\documentclass[12pt]{article}
\usepackage{amsthm,amsmath,amsfonts,amssymb,array, graphicx, epsfig}
\usepackage[latin1]{inputenc}

\newcommand{\be}{\begin{enumerate}}
\newcommand{\ee}{\end{enumerate}}
\newcommand{\beq}{\begin{eqnarray*}}
\newcommand{\eeq}{\end{eqnarray*}}
\newcommand{\beqnr}{\begin{eqnarray}}
\newcommand{\eeqnr}{\end{eqnarray}}
\newcommand{\llll}{\langle \langle}
\newcommand{\rrrr}{\rangle \rangle}
\newcommand{\rrr}{\rangle \rangle}

\newcommand{\beg}{\begin{equation}}
\newcommand{\ed}{\end{equation}}

\newcommand{\feddich}{\hfill $\Box$\\}

\newcommand{\R}{\mathbb{R}}
\newcommand{\N}{\mathbb{N}}
\DeclareMathOperator{\diag}{diag}
\DeclareMathOperator{\trace}{tr}

\DeclareMathOperator{\dimens}{dim}

\newcommand{\spn}{\operatorname{span}}

\newtheorem{assumption}{Assumption}
\newtheorem{proposition}{Proposition}

\newtheorem{definition}{Definitions}
\newtheorem{theorem}{Theorem}
\newtheorem{Lemma}{Lemma}

\newtheorem{corollary}{Corollary}

\def\e{\epsilon}

\def\T{{\relax\ifmmode I\!\!\hspace{-1pt}T\else$I\!\!\hspace{-1pt}T$\fi}}

\def\cO{{\cal O}}
\def\cA{{\cal A}}

\def\cF{{\cal F}}

\def\bI{{\bf I}}

\def\cD{{\cal D}}
\def\cI{{\cal I}}
\def\cJ{{\cal J}}

\def\cT{{\cal T}}

\def\cV{{\cal V}}
\def\cL{{\cal L}}
\def\cB{{\cal B}}
\def\cG{{\cal G}}

\def\cR{{\cal R}}

\def\bM{{\bf M}}

\def\ll{\langle}
\def\rr{\rangle}

\newcommand{\XX}{\mathcal{X}}
\newcommand{\beqn}{\begin{equation}}
\newcommand{\eeqn}{\end{equation}}

\def\endproof{\feddich}

\newcommand{\bx}{\mathbf{x}}

\usepackage[font=small,format=plain,labelfont=bf,up,textfont=it,up]{caption}

\makeatletter \@addtoreset{equation}{section} \makeatother

\newcommand\eref[1]{(\ref{#1})}

\def\int{\intop\limits}

\newcommand\norm[1]{\left\|#1\right\|}

\newcommand\bU{{\bf U}}

\newif\ifNZB

\makeatletter
\newcommand{\tr}{{\mathop{\operator@font T}\nolimits}}
\makeatother

\usepackage[usenames]{color}

\newcommand{\unred}{\color{black}}
\usepackage{algorithm}
\usepackage{algorithmic}

\newcommand{\skapro} [1] {\left\langle#1\right\rangle} 
\newcommand{\Lska} [1] {\skapro{#1}_{L_2}}  
\newcommand{\sska} [1] {\skapro { \skapro{ #1 } }}

\newcommand{\JJJ}{\ensuremath{{\cal J}}}

\begin{document}
\title{Direct minimization for calculating invariant subspaces
in density functional computations  of the  electronic structure
\thanks{This work was supported by the DFG SPP 1445: ``Modern and universal first-principles methods for many-electron systems in chemistry and physics'' and the EU NEST project BigDFT.}
}
\author{R. Schneider, T. Rohwedder, A. Neelov, J. Blauert}
\date{}
\maketitle

\begin{abstract}
In this article, we analyse three related preconditioned steepest descent algorithms, which are partially popular in Hartree-Fock and Kohn-Sham theory as well as invariant subspace computations, from the viewpoint of minimization of the corresponding functionals, constrained by orthogonality conditions. We exploit the geometry of the of the admissible manifold, i.e. the invariance with respect to unitary transformations, to reformulate the problem on the Grassmann manifold as the admissible set. We then prove asymptotical linear convergence of the algorithms under the condition that the Hessian of the corresponding Lagrangian is elliptic on the tangent space of the Grassmann manifold at the minimizer.
\end{abstract}\bigskip

\section{Introduction}

On the length-scale of atomistic or molecular systems, physics is
governed by the laws of quantum mechanics. A reliable
computation required in various fields in modern sciences and
technology should therefore be based on the first principles
of quantum mechanics, so that \textit{ab initio} computation of the
electronic wave function from the stationary electronic Schr\"odinger 
equation is a major working horse for many applications in this area. 
To reduce computational demands, the high dimensional problem of computing 
the wave function for $N$ electrons is often, for example in Hartree-Fock and Kohn-Sham theory,
replaced by a nonlinear system of equations for a set $\Phi = (\varphi_1, \ldots, \varphi_N)$
 of single particle wave functions $\varphi_i(\bx) \in V = H^1 ( \mathbb{R}^3)$. 
This ansatz corresponds to the following abstract formulation 
for the minimization of a suitable energy functional $\cJ$:

{\bf Problem 1}: Minimize
\beqnr \label{nonlinminiprob}
\mathcal{J}: V^N \to \mathbb{R} , ~~~~\mathcal{J}   (\Phi) ~~
= ~\mathcal{J}   (\varphi_1, \ldots, \varphi_N)  ~~\longrightarrow~~ \mbox{min},
\eeqnr
where $\cJ$ is a sufficiently often differentiable functional which is 
\be
\item[(i)] invariant with respect to
unitary transformations, i.e.
\beqnr
\mathcal{J} ( \Phi ) = \mathcal{J} ( \Phi \mathbf{U} ) = \mathcal{J}
( ( \sum_{j=1}^N u_{i,j} \phi_j)_{i=1}^N ) ,
\label{uniinvari}
\eeqnr
for any orthogonal matrix $\bU \in \R^{n\times n}$,
and 
\item[(ii)]subordinated to the orthogonality constraints
\beqnr
\langle \varphi_i , \varphi_j \rangle ~~:= ~~\int_{\mathbb{R}^3 } \varphi_i  ({x} )
 \varphi_j ({x} ) d{x} ~~=~~\delta_{i,j}.
\label{orthocon}
\eeqnr
\ee
In the present article, we shall be concerned with minimization techniques for $\cJ$ along the admissible
manifold characterized by \eref{orthocon}.  The first step towards this will be to set up the theoretical framework of the \emph{Grassmann manifold} to be introduced in section \ref{theoframesec}, reflecting the constraints (i) and (ii) imposed on the functional $\cJ$ and the minimizer $\Phi$, respectively.
In applications in electronic structure theory, formulation of the first order optimality (necessary) condition for the problem \eref{nonlinminiprob} results in a nonlinear eigenvalue problem of the kind:
\beqnr
A_{\Phi }  \varphi_i  = \lambda_i \varphi_i ,~~~~ \lambda_1
 \leq \lambda_2 \leq\ldots \leq \lambda_N \label{eigenveq}
\eeqnr
for $N$ eigenvalues $\lambda_i$ and the corresponding solution functions assembled in $\Phi$. 
In these equations, the operator $  A_{\Phi }$, is a symmetric bounded linear mapping $  A_{\Phi } : V= H^1 (
\mathbb{R}^3) \to V' = H^{-1} ( \mathbb{R}^3) $ depending on $\Phi$, so that we are in fact faced with a nonlinear eigenvalue problem. $  A_{\Phi }$ is called the \emph{Fock operator} in Hartree-Fock theory, and  \emph{Kohn-Sham} Hamiltonian in \emph{density
functional theory} (DFT) respectively. We will illustrate the relation between \eref{eigenveq} and the minimization task above in further detail in section \ref{schrosec}. In this work, our emphasis will rather be on the algorithmic
approximation of the minimizer of $\cJ$ , i.e. an invariant subspace $\spn [\Phi] := \spn\{ \varphi_1 , \ldots , \varphi_N \},$ of \eref{eigenveq}, in the corresponding energy space $V^N$ than on computation of the eigenvalues $\lambda_1, \ldots , \lambda_N$.\medskip

One possible procedure for computing the minimum of $\cJ$ is the so-called \textit{direct minimization}, utilized e.g.
in DFT calculation, which performs a steepest descent algorithm by updating the gradient of $\cJ$, i.e. the  Kohn-Sham
Hamiltonian or Fock operator, in each iteration step. Direct minimization, as proposed in
\cite{cg}, is prominent in DFT calculations if good preconditioners
are available and the systems under consideration are large, e.g. for the computation of
electronic structure in
 bulk crystals using
plane waves, finite differences \cite{becker} and the recent wavelet code
developed in the BigDFT project (see \cite{bigdft}). In contrast to the direct minimization procedure
 is the {\em self consistent field iteration (SCF)},  which keeps the Fock operator fixed until convergence of the
corresponding eigenfunctions and updates the Fock operator
 thereafter, see section \ref{schrosec}.


\smallskip

In the rest of this article, we will pursue different variants of projected gradient algorithms to be compiled in section \ref{algorithsec}. In addition, we will (for the case where the gradient $\cJ'(\Phi)$ can be written as an operator $A_{\Phi}$ applied to $\Phi$, as it is the case in electronic structure calculation) investigate an algorithm based on \cite{arias}
following a preconditioned  steepest descent along geodesics on the manifold.
so that no re-projections onto the admissible manifold are required. 
It turns out that all these algorithms to be proposed perform in a similar way.
For matters of rigorous mathematical analysis, let us note at this point that
  the mathematical theory about Hartree-Fock is
still too incomplete to prove the assumptions required in the present paper;
 even less is known for Kohn-Sham
equations, due to the fact that there are so many different models
used in practice. If the assumptions are not met for a particular
problem, it is not clear whether it is a deficiency of the problem or
a real pathological situation. Along with \eref{nonlinminiprob}, we will therefore 
consider the following simplified prototype problem for a \textit{fixed} operator $A$:\medskip

{\bf Simplified problem 2:} Minimize
\beqnr
\mathcal{J}_A   (\varphi_1, \ldots, \varphi_N) :=
\sum_{i=1}^N \langle\varphi_i ,  A \varphi_i \rangle~~\longrightarrow~~ \min, ~~~~
\langle \varphi_i , \varphi_j \rangle ~~= ~~\delta_{i,j}.\label{miniprob} \eeqnr

Analogous treatment with Lagrange techniques shows that this special case of problem 1 is the problem of computing the first $N$ eigenfunctions,
resp. the lowest $N$ eigenvalues of $A$ (see Lemma \ref{prop:localminima1}).
While this is an interesting problem by itself, e.g. if $ \lambda $ is an eigenvalue of multiplicity $N$,
 it is also of interest as a sort of prototype: Properties that can
be proven for this problem may hold in the more general
case for Hartree-Fock or Kohn-Sham. In particular, we will show that for $A$ symmetric and bounded from below, the Hessian of the Lagrangian, taken at the solution $\Psi$, is elliptic on a specific tangent manifold at $\Psi$, an essential ingredient to prove linear convergence of all of the proposed algorithms in section \ref{convsec}. The same convergence results will be shown to hold for \eref{nonlinminiprob} if we impose this ellipticity condition on the Lagrangian of $\cJ$ of the nonlinear problem.
Note that the problem type \eref{miniprob} also arises in many other
circumstances, which we will not consider here in detail. Let us just note that the algorithms presented in section \ref{algorithsec} also provide reasonable routines for the inner cycles of the SCF procedure.\\
In the context of eigenvalue computations, variants of our basic algorithm 1, applied to problem 2,
have been considered by several authors (see e.g.
\cite{knyazev, knyaneymeyr, neymeyr}) reporting excellent performance, in
particular if subspace acceleration techniques are applied and the
preconditioner is chosen appropriately; in \cite{rsz, drsz}, an adaptive variant was recently proposed and analysed for the simpler case $N=1$.  In contrast to all these papers, we will view the algorithms 
as steepest descent algorithms for
optimization of $\cJ$ under the orthogonality constraints given above, as such a
systematic treatment does not only simplify the proofs but also 
provides the insight necessary to understand
the direct minimization techniques for the more complicated
nonlinear problems of the kind \eref{nonlinminiprob} in DFT and HF. \\
Our analysis will cover closed (usually finite dimensional) subspaces of $V_h \subset V$ as well 
as the energy space $V$ itself, so that finite dimensional approximations by Ritz-Galerkin methods and also
finite difference approximations are included in our analysis. In particular, our results are also valid if
Gaussian type basis functions are used. The convergence 
rates will be independent of the discretization parameters like mesh
size.  However, the choice of an appropriate preconditioning mapping to be used in our algorithms is crucial.
Fortunately, such preconditioners can often easily be constructed, e.g. by the use of multigrid methods for
finite elements, finite differences or wavelets, polynomials \cite{hackbusch, cg, becker}. 
Our analysis will show that for the gradient algorithms under consideration, it suffices to use a \emph{fixed} preconditioner respectively relaxation parameter. In particular, no expensive line search is required.\\
All results proven will be local in nature meaning that the
initial guess is supposed to be already  sufficiently close to the
exact one. At the present stage, we will for the sake of simplicity consider only real valued solutions for the minimization problem. Nevertheless, complex
valued functions can be treated by minor modifications. Note that since the present approach is completely based on a variational
framework, i.e. considering a constrained optimization problem, it does not include unsymmetric eigenvalue problems or the computation of other eigenvalues than the lowest ones.

\section{Optimization on Grassmann manifolds}\label{theoframesec} 

The invariance of the functional $\cJ$ with respect to uniform transformations among the eigenfunctions shows a certain redundance inherent in the formulation of the minimization task \eref{nonlinminiprob}. Therefore, it will be more advantageous to factor out the unitary invariance of the functional $\cJ$, resulting
in the usage of the Stiefel and Grassmann manifolds, originally defined in finite dimensional Euclidean Hilbert
spaces in \cite{arias}, see also \cite{matrixman} for an extensive exposition. In this section, we will generalize this concept for the present infinite dimensional space $V^N$ equipped with the $L_2$ inner product. In the next section, we will then apply this framework to the minimization problems for the HF and KS functionals. First of all, we shall briefly introduce the spaces under consideration and some notations.

\subsection{Basic notations} 

Letting $H= L_2 := L_2 ( \mathbb{R}^3) $ or a closed subspace of $L_2$, we will work with a Gelfand triple $ V \subset H \subset V'$ with the usual $L_2$ inner product $\ll. , ..\rr,$ as dual pairing on $ V' \times V$, where either $V := H^1 = H^1(\R^3)$ or an appropriate subspace corresponding to a Galerkin discretization. Because the ground state is determined by a set $\Phi$ of $N$ 
one-particle functions $\varphi_i \in V$, we will formulate the optimization problem on an admissible subset of $V^N$.
To this end, we extend inner products and operators from $V$ to $V^N$ by the following

\begin{definition} 
For $\Psi = (\psi_1 , \ldots , \psi_N)\in V^N,$ $\Phi = (\varphi_1 , \ldots , \varphi_N) \in (V^N)'=(V')^N$, and the $L_2$ inner product $\ll.,.. \rr$ given on $H = L_2$,
we denote 
$$\skapro{\Phi^T \Psi} := ( \skapro{\varphi_i, \psi_j} )_{i,j=1}^N \in \R^{N \times N} \ ,$$
and introduce the dual pairing 
$$\sska { \Phi, \Psi } := \trace  \skapro{\Phi^T \Psi}   = \sum_{i=1}^N
\skapro{ \varphi_i, \psi_i}$$  on  $  (V')^N\times V^N$.

Because there holds $V^N = V \otimes \R^N$, we can canonically expand any operator $R: V \to V'$ to
an operator 
\beqnr \cR := R \otimes I: V^N = V \otimes \R^N \to V'^N, \Phi 
\mapsto \cR \Phi = (R \varphi_1, \ldots, R \varphi_N). \label{calliVN} \eeqnr
Throughout this paper, for an operator $V \to V'$ denoted by a capital letter as $A, B, D, \ldots$, the same calligraphic letter $\cA, \cB, \cD, \ldots$, will denote this expansion to $V^N$.\end{definition} 

Further, we will make use of the following operations:
\begin{definition} 
 For $\Phi \in V^N$ and $ \bM \in \R^{N\times N}$, we define 
the set $\Phi \bM= (I \otimes \bM)\Phi$ by
$(\Phi \bM)_j := \sum_{i=1}^N  m_{i,j} \varphi_i$, cf. also the notation in \eref{uniinvari},
and for $\phi \in  V$ and $v=(v_1, \ldots, v_N) \in \R^N$ the element $\phi \otimes v \in V^N$ by $(v_1\phi, \ldots, v_N\phi).$
Finally, we denote by $O(N)$ the orthogonal group of $\R^{N \times N}$.
\end{definition}

\subsection{Geometry of Stiefel and Grassmann manifolds}\label{grasstiefel}

Let us now introduce the admissible manifold and prove some of its basic properties.
Note in this context that well established results of \cite{arias} for the case in the finite dimensional Euclidean spaces cannot
be applied to our setting without further difficulties, because the norm induced by
the $ L_2 $ inner product is weaker than the present $V$-norm.\\

Our aim is to minimize the functionals $\cJ(\Phi)$, where $\cJ$ is either $\cJ_{HF}$, $\cJ_{KS}$ or $\cJ_A$, 
under the
orthogonality  constraint $\langle \varphi_i , \varphi_j \rangle = \delta_{i,j}$, i.e.
\begin{eqnarray} \label{ortho_constraint}
\skapro{\Phi^T \Phi} = {\bf I} \in \R^{N \times N}.
\end{eqnarray}
The subset of $V^N$  satisfying the property (\ref{ortho_constraint}) 
is called  the \emph{Stiefel manifold} (cf. \cite{arias})
$$\mathcal{V}_{V,N} := \{ \Phi
=(\varphi_i)_{i=1}^N | \varphi_i \in V, \skapro{\Phi^T \Phi}  - \bI = \mathbf{0}
\in \R^{N \times N} \} \, ,$$ 
i.e. the set of all
orthonormal bases of $N$-dimensional subspaces of $V$.

All functionals $ \mathcal{J}$ under consideration
are unitarily invariant, i.e. there holds \eref{uniinvari}. To get rid of
this nonuniqueness, we will identify all orthonormal bases $ \Phi \in
\mathcal{V}_{V,N}$ spanning the same subspace $ V_{\Phi}
:= \mbox{span } \{ \varphi_i : i=1, \ldots , N \} $. To this end we
consider the \emph{Grassmann manifold}, defined as the quotient
$$\mathcal{G}_{V,N} := \mathcal{V}_{V,N} / {\sim}$$
 of the Stiefel
manifold with respect to the equivalence relation $\Phi {\sim}
\widetilde{\Phi}$ if $\widetilde{\Phi} = \Phi \bU$ for any $\bU \in O(N)$. We usually omit the indices and write $\cV$ for
$\mathcal{V}_{V,N}$, $\cG$ for $\mathcal{G}_{V,N}$ respectively. To
simplify notations we will often also work with representatives instead of
equivalence classes $[\Phi]\in \cG$.

The interpretation of the Grassmann manifold as equivalence classes of
orthonormal bases spanning the same $N$-dimensional subspace is just one way
to define the Grassmann manifold. We can as well identify the subspaces with
orthogonal projectors onto these spaces. To this end, let us for $\Phi =(\varphi_1, \ldots, \varphi_N) \in \cV^N$ denote by $D_{\Phi}$
the $L_2$-orthogonal projector onto $\spn \{ \varphi_1, \ldots, \varphi_N \}$. It is straightforward to verify

\begin{Lemma}
There is a one to one relation identifying $\cG $ with the set of rank $N$
$L_2$-orthogonal projection operators $ D_{\Phi} $.
\end{Lemma}

In the following, we will compute the tangent spaces of the manifolds defined above for later usage.

\begin{proposition}\label{stiefeltangentspace}
The tangent space of the Stiefel manifold at $\Phi \in \cV$ is given by
$$ \mathcal{T}_{\Phi}\mathcal{V} = \{ X \in V^N \, | \, \skapro{X^T \Phi} = - \skapro{\Phi^T X} \in \R^{N \times N} \} \, .$$
The tangent space of the Grassmann
manifold is
 \begin{eqnarray*}
\mathcal{T}_{[\Phi]} \mathcal{G} & = &
 \{ W \in V^{N}| \skapro{W^T \Phi} = \mathbf{0} \in \R^{N \times N} \} \\
& = & (\spn\{\varphi_1, \ldots, \varphi_N\}^{\perp})^N \, .
\end{eqnarray*}
Thus, the operator $( \cI - \cD_{\Phi})$, where $D_\Phi$ is the $L_2$-projector onto the space spanned by $\Phi$, is an  $L_2 $-orthogonal
projection from $V^N$ onto the tangent space $\cT_{[\Phi]} \mathcal{G}$.
\end{proposition}
\begin{proof}
If we compute the Fr\'{e}chet derivative of the constraining condition 
$$g(\Phi) := \skapro{ \Phi^T \Phi} - \bI = \mathbf{0}$$ for the Stiefel manifold,
the first result follows immediately.
To prove the second result, we consider the quotient structure of the Grassmann manifold 
and decompose the tangent space $\mathcal{T}_{\Phi}  \mathcal{V}$ of the  Stiefel 
manifold at the representative $\Phi$ into a component tangent to the
set $[\Phi]$, which we call the \emph{vertical space}, and a component containing the
elements of $\mathcal{T}_{\Phi}  \mathcal{V}$ that are orthogonal to the vertical space, 
the so-called \emph{horizontal space}.
If we move on a curve in the Stiefel manifold with direction in the vertical space, we do
not leave the equivalence class $[\Phi]$. Thus only the horizontal space defines the 
tangent space of the quotient
$\mathcal{G} = \mathcal{V} / O(N) $. 
The horizontal space is computed in the following lemma, from which the claim follows.
\end{proof}

\begin{Lemma}
The vertical space at a point $\Phi \in \mathcal{V}$  (introduced in the proof of proposition \ref{stiefeltangentspace}) is the 
set 
$$\{ \Phi \bM | \bM = -\bM^T \in \mathbb{R}^{N \times N} \} \, .$$
The horizontal space is given by
$$\{ W \in V^{N}| \skapro{W^T \Phi} =  \mathbf{0} \in \R^{N \times N} \} \ .
$$
\end{Lemma}
\begin{proof}
To compute the tangent vectors of the set $[\Phi]$, we consider 
a curve $c(t)$ in $[\Phi]$ emanating from $\Phi$. Then $c$ is of the form 
$c(t) = \Phi \bU(t)$ for a curve $\bU(t) \in O(N)$ with $\bU(0)= \bI_{N \times N}$. 
Differentiating $\bI_{N \times N} = \bU(t) \bU(t)^T$ at $t=0$ yields 
$\bU'(0) = - \bU'(0)^T$ and we get that every vector of the vertical space 
is of the form $\Phi \bM$ where $\bM$ is skew symmetric.\\
Reversely, for any skew symmetric matrix $\bM$ we find a curve $\bU(t)$ in $O(N)$ 
emanating from $\Phi$ with direction $\bM$, and $c(t) := \Phi \bU(t)$ is a curve 
with direction $\dot{c}(0) = \Phi \bM$, and thus the first assertion follows.

To compute the horizontal space, we decompose $W \in \mathcal{T}_{\Phi} \mathcal{V}$ into
$W=\Phi \bM + W_{\perp}$, where $W_{\perp} := W  - \Phi \skapro{\Phi^T W} \in \Phi^{\perp}$, 
$\bM := \skapro{\Phi^T W}$.
Then $\bM$ is  an antisymmetric matrix, which implies that $\Phi \bM$ is in the vertical space,
and that  the horizontal space is given by all $\{W_{\perp} = W - \Phi \skapro{\Phi^T W} | W \in \mathcal{T}_{\Phi} \mathcal{V} \}$.
Let us note that this set is the range of the operator $(I-\cD_{\Phi})$. This operator is continuous and of finite codimension.
If $W_{\perp}= W - \Phi \skapro{\Phi^T W} $ is in the horizontal space, then
$$\skapro{W_{\perp}^T \Phi} = \skapro{W^T \Phi} - \skapro{\Phi^T \Phi}  \skapro{W^T \Phi} = \mathbf{0}.$$

Reversely, if $W \in V^N$ with $\skapro{W^T \Phi}=\mathbf{0}$, then $W$
is in $\cT_{\Phi} \cV$ and from
$(\cI-\cD_{\Phi})W = W - \Phi \skapro{\Phi^TW} = W$ we get that $W$ is in
the range of $I-\cD_{\Phi}$, being the  $L_2 $-orthogonal
projection from $V^N$ onto the tangent space $T_{[\Phi]} \mathcal{G}$.
\end{proof}

To end this section, let us prove a geometric result needed later.

\begin{Lemma}\label{lem:projdiffer}
Let $[\Psi] \in \cG$, $D^{*}$ the $L_2$-projector on $\spn [\Psi]$, $\cD^{*}$ is its expansion as above and $||.||$ is the norm induced by the $L_2$ inner product. For any $\Phi = (\varphi_1, \ldots, \varphi_N ) \in \cV$ sufficiently close to $[\Psi] \in \cG$ in the sense that for all $i \in \{1,\ldots,N\}$, $||(I -D^{*})\varphi_i||<\delta$, there exists an orthonormal basis $\bar{\Psi}  \in \cV$  of $\spn [\Psi]$ for which
$$ \Phi - \bar{\Psi} ~~=~~ (I- \cD^{*})\Phi  ~+~ \cO(||(I- \cD^{*})\Phi||^2).$$ 
\end{Lemma}
\begin{proof}
For $i = 1, \ldots, N$, let
$$ \widetilde{\psi_i} ~~=~~ \arg\min\{ ||\psi -\varphi_i||, \psi \in \spn\{\psi_i|i=1,\ldots ,N\}, || \psi|| = 1\} ~~=~~ D^{*}\varphi_i/||D^{*}\varphi_i||,$$ and set $\widetilde{\Psi} := (\widetilde{\psi}_1, \ldots, \widetilde{\psi}_N).$ If we denote by $\widetilde{P}_i$ the $L_2$ projector on the space spanned by $\widetilde{\psi}_i$, it is straightforward to see from the series expansion of the cosine that 
\beqnr (I- D^{*})\varphi_i ~~=~~ (I- \widetilde{P}_i)\varphi_i ~~=~~ \varphi_i - \widetilde{\psi_i}~~ +~~ \cO(||(I- D^{*})\varphi_i||^2) \label{quadsing}\eeqnr
The fact that $\widetilde{\Psi} \notin \cV$ is remedied by orthonormalization of $\widetilde{\Psi}$ by the Gram-Schmidt procedure. For the inner products occurring in the orthogonalization process (for which $i\neq j$), there holds
\beq \ll\widetilde{\psi}_i,\widetilde{\psi}_j \rr ~&=&~ \ll\widetilde{\psi}_i - \varphi_i,\widetilde{\psi}_j \rr + \ll \varphi_i,\widetilde{\psi}_j - \varphi_j \rr
 + \ll \varphi_i,  \varphi_j\rr\\
~&=&~ - \ll(I- D^{*})\varphi_i,\widetilde{\psi_j}\rr -  \ll(I- D^{*})\varphi_i , (I- D^{*})\varphi_j \rr ~ +~ \cO(||(I- D^{*})\varphi_i||^2). \\
~&=&~  \cO(||(I- \cD^{*})\Phi||^2)  \eeq
where we have twice replaced $\varphi_i - \widetilde{\psi_i}$ by $(I- D^*)\varphi_i$ according to \eref{quadsing} and made use of the orthogonality of $D^*$. In particular, for $\Phi$ sufficiently close to $[\Psi]$, the Gramian matrix is non-singular because the diagonal elements converge quadratically to one while the off-diagonal elements converge quadratically to zero.  By an easy induction for the orthogonalization process and a Taylor expansion for the normalization process,
we obtain that $\widetilde{\Psi} $ differs from the orthonormalized set $\bar{\Psi}:= (\bar{\psi_1}, \ldots, \bar{\psi_N})$ only by a error term depending on $||(I- \cD^{*})\Phi||^2$.
 Therefore, 
\beq \varphi_i  - \bar{\psi}_i ~&=&~~ \varphi_i - \widetilde{\psi_i}~ +~ \cO(||(I- \cD^{*})\Phi||^2)~~=~~ (I- D^{*})\varphi_i ~ +~ \cO(||(I- \cD^{*})\Phi||^2), \eeq
so that  \beq \Phi - \bar{\Psi} ~~=~~ (I- \cD^{*})\Phi ~ +~ \cO(||(I- \cD^{*})\Phi||^2), \eeq
and the result is proven.
\end{proof}

\subsection{Optimality conditions on the Stiefel manifold}\label{ministiefel}

By the first order optimality condition for minimization tasks, a minimizer $ [\Psi]  \in\cG$
of the functional
 $\cJ: \cG \to \R, \Phi \mapsto \cJ ( \Phi ) $ over the Grassmann manifold
$\cG$ satisfies 
\begin{equation}\label{eq:1ordercond}
\langle \langle \cJ' (\Psi ) , \delta \Phi \rangle \rangle=
0 ~~~ \makebox{for all}~~ \delta \Phi \in \mathcal{T}_{[\Psi]} \cG \ ,
\end{equation}
i.e. the gradient  $\cJ' ( \Psi ) \in (V')^N= (V^N)'$ vanishes on the
tangent space $ \mathcal{T}_{\Psi} \cG$ of the Grassmann manifold.
This property can also be formulated by
$$
 \skapro{(\delta \Phi)^T   \cJ' (\Psi )} = \mathbf{0} ~~~~ \makebox{for all}~~ \delta \Phi \in
\mathcal{T}_{[\Psi]} \cG,
$$
or equivalently, by Lemma \ref{stiefeltangentspace},
\beqnr
 \llll(\cI - \cD_{\Psi}) \cJ' (\Psi ), \Phi \rrr = 0 ~~~~ \makebox{for all}~~ \Phi \in V^N,
\label{firstoptimality}
\eeqnr
that is, in strong formulation,
\beqnr(\cI - \cD_{\Psi})\cJ' (\Psi )~=~ \cJ' (\Psi ) - \Psi \Lambda  ~=~ 0 \in (V')^N, \label{strongopt}\eeqnr where $\Lambda = (\ll(\cJ'(\Psi))_j, \psi_i\rr)_{i,j=1}^N$
and $(\cJ'(\Psi))_i \in V'$ is the $i$-th component of $\cJ'(\Psi)$. Note that this corresponds to one of the optimality conditions for the Lagrangian yielded from the common approach of the Euler-Lagrange minimization formalism: Introducing the Lagrangian
\beqnr \cL(\Phi,  \Lambda )
:= \frac{1}{2} \left( \cJ( \Phi ) + \sum \lambda_{i,j}
(\Lska{\varphi_i, \varphi_j}- \delta_{i,j} ) \right),\eeqnr
the condition for the derivative restricted to $V^N$, here denoted by $\mathcal{L}^{(1,\Psi)}(\Psi, \Lambda)$ for convenience, is given by
\beqnr \mathcal{L}^{(1,\Psi)}(\Psi, \Lambda) &=&  \mathcal{J}' ( \Psi ) - (\sum_{k=1}^{N}
\lambda_{i,k} \psi_k)_{i=1}^N = 0 \in (V')^N. \label{lagrangian}\eeqnr
Testing this equation with $\psi_j, j=1, \ldots, N$, verifies the Lagrange multipliers indeed agree with the $\Lambda$ defined above, so that \eref{firstoptimality} and \eref{lagrangian} are equivalent.
Note also that the remaining optimality conditions, \beq \frac{\partial{\cL}}{\partial \lambda_{i,j}} &=&
\frac{1}{2} \left( (\Lska{\psi_i, \psi_j}- \delta_{i,j} \right)~=~0,\eeq of the Lagrange formalism are now incorporated in the framework of the Stiefel manifold.\\
From the representation \eref{strongopt}, it follows that the Hessian $\cL^{(2, \Psi)}(\Psi, \Lambda)$ of the Lagrangian \eref{lagrangian}, taken at the minimum $\Psi$ and with the derivatives taken with respect to $\Psi$, is given by
$$\cL^{(2, \Psi)}(\Psi, \Lambda) \Phi ~=~ \cJ'' (\Psi )\Phi - \Phi \Lambda.$$
As a necessary second order condition for a minimum, $\cL(\Psi, \Lambda)^{(2, \Psi)}$ has to be positive semidefinite on $\cT_{[\Psi]}\cG$. 
For our convergence analysis, we will have to impose the stronger condition on $\cL^{(2, \Psi)}(\Psi, \Lambda)$ being elliptic on the tangent space, i.e.
\beqnr
\langle\langle \cL^{(2, \Psi)} ( \Psi, \Lambda ) \delta \Phi ~,~ \delta \Phi   \rangle\rangle ~~\geq~~
\gamma ~\| \delta \Phi \|_{V^N}^2, ~~~~ \makebox{for all}~~~ \delta \Phi \in \mathcal{T}_{[\Psi]}\cG\label{ellichap2}.
\eeqnr
It is an unsolved problem if this condition holds in general for the minimization problems of the kind \eref{nonlinminiprob} or if it depends on the functional under consideration; in particular, it is not clear whether it holds for the functionals of Hartree-Fock and density functional theory. In the case of Hartree-Fock, it suffices to demand that $\cL^{(2,\Psi)}( \Psi , \Lambda) > 0$ on $\cT_{[\Psi]}\cG$ because this already implies $\cL^{(2,\Psi)}( \Psi , \Lambda)$ is bounded away from zero, cf \cite{maday}. For the simplified problem, we will show in Lemma \ref{linearprobelli} that the assumption holds for symmetric operators $A$ fulfilling a certain gap condition.

\section{Minimization tasks in electronic structure\\ calculations}\label{schrosec}

We will now particularize the results of the last section to the functionals common in electronic structure calculation.
As the following section will show, the applications of interest in electronic structure calculations deal with the minimization of functionals $\cJ$ for which the gradient can be written as $\cJ'(\Phi) = \cA_{\Phi}\Phi,$ where $A_{\Phi}: V \to V'$ (and $\cA_{\Phi}$ its extention to $V^N$ by \eref{calliVN}). We conjecture that if the functional $\cJ$ only depends on the electronic density, that is, if condition \eref{uniinvari} holds, 
this form of $\cJ(\Phi)$ is valid in general, i.e. for each $\Phi \in \cG,$ there is an operator $A_{\Phi}$ so that $\cJ' (\Phi) = \cA_{\Phi}\Phi$. Nevertheless, we decided to formulate the algorithms (except algorithm 3) for $\cJ' (\Phi)$ rather than for $A_{\Phi}$ to emphasize the minimization viewpoint we pursue in this work and to display that the concrete structure of the Fock or Kohn-Sham operators does not enter anywhere in the proof of convergence given in section \ref{convsec}.\\
In this section, we will remind the reader of some basic facts about Hartree-Fock and Kohn-Sham theory, where our emphasis will be on the ansatzes leading to the problem of minimizing a nonlinear functional \eref{nonlinminiprob}. Also, we will review the concrete form the operator $\cJ'(\Phi) = \cA_{\Phi}\Phi$ in \eref{eigenveq} has in these applications.
For a more detailed introduction to electronic structure calculations, we refer the reader to the standard literature \cite{lebris, lebris4, helgaker, szabo}.
At the end of this section, we will investigate the simplified problem \eref{miniprob} and its connection to eigenvalue computations.

\subsection {Hartree-Fock and Kohn-Sham energy functionals\\ in quantum chemistry}\label{hfsub}

The commonly accepted model to describe atoms and molecules is by
means of the Schr\"odinger equation, which is in good agreement with
experiments as long as the energies remain on a level at which
relativistic effects can be neglected. We are mainly interested
in the stationary ground state of quantum mechanical systems,
given by the eigenfunction belonging to the lowest eigenvalue of the Hamiltonian $H$ of the system.
In the Born-Oppenheimer approximation
the Hamiltonian of
the (time-independent) {\em electronic Schr\"odinger equation} $H\Psi = E \Psi$ is 
given by
$$ {H} := - \frac{1}{2} \sum_{i=1}^{N^*} \Delta_i -
\sum_{i=1}^{N^{*}} \sum_{\nu = 1}^M \frac{Z_{\nu}}{\norm{x_i - R_{\nu} }}
 + \frac{1}{2} \sum_{i,j = 1, i\ne j}^{N^{*}} \frac{1}{ \norm{ x_i -x_j }} \, .$$
Here, $N^{*}$ denotes the number of electrons, $M$ the number of the nuclei, and $Z_{\nu}$, $R_{\nu}$ the charge respectively the coordinates
of the nuclei, which are the only fixed input parameters of the system. Note that we use atomic units, so that no physical constants appear in the Schr\"odinger equation. We also
neglect the interaction energy between the nuclei, since for a given constellation $(R_1, \ldots, R_M)$ of the
$M$ nuclei this
only adds a constant to the energy eigenvalues. Due to the Pauli principle for fermions, the wave function is required to be antisymmetric with respect to permutation of particle coordinates.
It is easy to see that every such antisymmetric solution
can be represented by a convergent sum of Slater determinants of the form
$$\psi_{SL}^{\Phi}(x_1, s_1  \ldots, x_{N^{*}}, s_{N^{*}}) := \frac{1}{\sqrt{N^{*}!}}
\det(\varphi_i(x_j,s_j )),~~x_i \in \R^3, ~~s_i = \pm \frac 12$$ where $\Phi = (\varphi_i)_{i=1}^{N^{*}} \in H^1(\R^3 \times \{\pm \frac 12\})^{N^{*}}$ and $\skapro {\varphi_i,
\varphi_j} = \delta_{i,j}$. In \emph{Hartree-Fock (HF) theory}, one approximates the ground state of the system by
 minimizing the Hartree-Fock energy functional $\Phi \mapsto \JJJ_{HF}(\Phi) :=\skapro{
H \psi_{SL}^{\Phi},\psi_{SL}^{\Phi} }$ 
over the set of all wave functions consisting of \textit{one single} Slater
determinant $\psi_{SL}^{\Phi}(x_1, s_1  \ldots, x_{N^{*}}, s_{N^{*}}).$
Additional simplification is made by the {\em Closed Shell Restricted Hartree-Fock model} (RHF), 
given in a spin-free formulation for $N = N^*/2$ \emph{pairs} of electrons, so that $\Phi = (\varphi_i)_{i=1}^N \in H^1(\R^3)^N =: V^N$. Abbreviating $V(x) := -\sum_{\nu = 1}^M \frac{Z_{\nu}}{\norm{x - R_{\nu} }},$ the corresponding functional reads
\begin{eqnarray}
\cJ_{HF} ( \Phi )  :=   \sum_{i=1}^N 
\int_{\R^3} &\Big(&\frac{1}{2} |\nabla \varphi_i (x)|^2 + V(x) |\varphi_i (x)|^2
     ~~+~~ \frac{1}{2}\sum_{j=1}^N \int_{\R^3}
     \frac{|\varphi_j(y)|^2} { \norm{ x-y}} \, dy \, |\varphi_i(x)|^2 \nonumber \\
     &-&\frac{1}{2}\sum_{j=1}^N
     \int_{\R^3} \frac{\varphi_i(x) \varphi_j(x) \varphi_j(y) \varphi_i(y) }{\norm{x-y}} \, dy \Big) ~~ dx.\label{hffunc}
\end{eqnarray}

A minimizer of $\JJJ_{HF}$ is named Hartree-Fock ground state. Its existence has been proven in the case that
$\sum_{\mu=1}^K Z_{\mu}> N-1$ (\cite{liebsimon}, \cite{lions}). \\

The energy functional of the \emph{Kohn-Sham (KS) model} can be derived from the Hartree-Fock energy functional by two modifications: First of all, as a consequence of the Hohenberg-Kohn theorem (cf. \cite{dft}), it is formulated in terms of the electron density $n(x) = \sum_{i=1}^N |\varphi_i(x)|^2$ rather than in terms of the single particle functions; secondly, it replaces the nonlocal and therefore computationally costly exchange term in the Hartree-Fock functional (i.e. the fourth term in \eref{hffunc}) by an additional (a priori unknown) exchange correlation energy term $E_{xc}(n)$ also depending only on the electron density. The resulting energy functional reads 
$$\cJ_{KS} (\Phi )  =  \frac{1}{2} \sum_{i=1}^{N} \int_{\R^3} |\nabla \varphi_i(x)|^2 dx
            + \int_{\R^3} n(x) V(x) + \frac{1}{2} \int_{\R^3} \int_{\R^3} \frac{n(x) n(y) }{\norm{x-y}} \, dx  dy
            + E_{xc}(n). $$
Determining the ground state energy of Kohn-Sham theory then consists in a minimization of $\mathcal{J}_{KS}$ over all $\Phi = (\varphi_1, \ldots, \varphi_N) \in V^N$ with $\skapro{\varphi_i, \varphi_j} = \delta_{i,j}.$
Since the exchange correlation energy $E_{xc}$ is not known explicitly,
further approximations are necessary. The most simple approximation
for $E_{xc}$
 is the  local density approximation (LDA, cf. \cite{dreizler}) defined as
$E_{xc}^{LDA} (n)= \int_{\R^3} n(x) \e_{xc}^{LDA} (n(x)) \, dx$,
where $\e_{xc}^{LDA}$
 denotes the exchange-correlation energy of a particle in an electron gas with density $n$. If we split this expression in an exchange and a correlation part, we get 
$$E_{xc}^{LDA} (n)~ ~=~~ E_{x}^{LDA} (n) + E_{c}^{LDA} (n) ~~=~~ \int_{\R^3} n(x) \e_{x}^{LDA} (n(x)) \, dx ~~+ ~~\int_{\R^3} n(x)
\e_{c}^{LDA} (n(x)) dx,$$
where in the exchange part, $\e_{x}^{LDA} (n) = - C_D n^{\frac{1}{3}}$ and $C_D := \frac{3}{4} (\frac{3}{\pi})^{1/3}$ is the Dirac constant.\\
For the correlation part  $E_{c}^{LDA} (n) $, the expression $\e_{c}^{LDA}$ is
analytically unknown, but can be calibrated e.g. by Monte-Carlo methods.
We note that a combination of both HF and density functional models, namely the hybrid B3LYP,
is experienced to provide the best results in benchmark computations.

\subsection{Canonical Hartree-Fock and Kohn-Sham equations}

For the HF and KS functionals, we can compute the derivative of $\cJ$ and the Lagrange multipliers at a minimizer explicitly.

\begin{proposition}\label{hfksoperators}
For the functional $\cJ_{HF}$ of Hartree-Fock, $\cJ_{HF}' ( \Phi ) = \mathcal{A}_{\Phi} \Phi \in (V')^N $, where $A_{\Phi} =  F^{HF}_{\Phi}: H^1(\R^3) \to H^{-1}(\R^3)$ is the so-called Fock operator and $\cA_{\Phi}$ is defined by $A_{\Phi}$ through \eref{calliVN}; using the notation of the \emph{density matrix}
\begin{eqnarray*}
\rho_{\Phi}(x,y) & :=& N \int_{\R^{3(N-1)} }\psi_{SL}^{\Phi}(x,x_2, \ldots, x_N) ~~ \psi_{SL}^{\Phi}(y,x_2, \ldots, x_N)  ~~ dx_2 \cdots dx_N 
\\ &=&  \sum_{i=1}^N \varphi_i(x) \varphi_i(y) 
\end{eqnarray*} and the electron density
 $n_{\Phi}(x):= \rho_{\Phi}(x,x)$ already introduced above.
It is given by
\begin{eqnarray*}
    F^{HF}_{\Phi} \varphi(x)  :=  
    & - & \frac{1}{2} \Delta \varphi (x) + V(x) \varphi(x) \\
    & + & \int_{\R^3} \frac{n_{\Phi}(y)} { \norm{ x-y}} \, dy \, \varphi(x)
    -  \int_{\R^3} \frac{ \rho_{\Phi}(x,y) \varphi(y) }{\norm{x-y}} \, dy.
\end{eqnarray*}
For the gradient of the Kohn-Sham functional $\cJ_{KS}$, there holds the following: Assuming that $E_{xc}$ in $\cJ_{KS}$ is differentiable and denoting by $v_{xc}$ the derivation of $E_{xc}$ with respect to the density $n$, we have $\cJ' ( \Phi ) = \mathcal{A}_{\Phi} \Phi \in (V')^N $, with $A_{\Phi} = F_{n}^{KS}$ the Kohn-Sham Hamiltonian, given by
$$F^{KS}_n \varphi_i := - \frac{1}{2} \Delta \varphi_i + V(x) \varphi_i + \left( n  \star \frac{1}{ \norm{ \cdot }} \right)  \varphi_i + v_{xc}(n) \varphi_i.$$\\
In both cases, the
Lagrange multiplier $\Lambda$ of \eref{lagrangian} at a minimizer $\Psi = (\psi_1, \ldots, \psi_N) $ is given by
\begin{eqnarray} \label{lagrange_multiplier}
\lambda_{i,j} = \ll A_{\Psi} \psi_i, \psi_j \rr.
\end{eqnarray}
There exists a unitary transformation $\bU \in O(N)$ amongst the functions $\psi_i,$ $i = 1, \ldots, N$ 
such that the Lagrange
multiplier is diagonal for $\Psi\bU = (\tilde{\psi_1}, \ldots, \tilde{\psi_N})$, 
$$\lambda_{i,j} := \ll A \tilde{\psi}_i, \tilde{\psi}_j \rr = \lambda_i \delta_{i,j}.$$
so that the ground state  of the HS resp. KS functional (i.e. minimizer of $\cJ$) satisfies the nonlinear \emph{Hartree-Fock} resp. \emph{Kohn-Sham} eigenvalue equations
\beqnr F^{HF}_{\Psi}\psi_i = \lambda_i \psi_i, ~~~\makebox{resp.}~~~F^{KS}_n \psi_i = \lambda_i \psi_i,~~~ ~~~\lambda_i \in \R, ~~~~~i=1, \ldots, N , \label{fockev}\eeqnr
for some  $\lambda_1, \ldots , \lambda_N \in \R$ and a corresponding set of orthonormalized functions $\Psi= (\psi_i)_{i=1}^{N}$ 
up to a unitary transformation $\bU$.
\end{proposition}

The converse result, i.e. if for a collection $\Phi = (\varphi_1, \ldots, \varphi_N)$ belonging to the $N$ lowest eigenvalues of the Fock operator in \eref{fockev}, the corresponding Slater determinant actually gives the Hartree-Fock energy by $\cJ(\Phi) = \ll H \psi_{SL}^{\Phi}, \psi_{SL}^{\Phi} \rr $, 
is not
known yet.

\subsection{Simplified problem}

The practical significance of the simplified problem \eref{miniprob} is given by the following result, which shows that for symmetric $A$, the  minimization of $\cJ_A$ is indeed equivalent to finding an orthonormal basis $\{ \psi_i : 1\leq i \leq N\}$ spanning the invariant subspace of $A$ given by the first eigenfunctions of $A$.

\begin{proposition}\label{prop:localminima1} Let $A$ in the simplified problem \eref{miniprob} a bounded symmetric
operator. The gradient of the functional $\cJ_A$ is then given by $\cJ' ( \Phi ) = \mathcal{A} \Phi \in (V')^N $.
Therefore, $\Psi$  is a stationary point of
$ \cL$ if and only if there exists an orthogonal transformation $\bU$ such that $\Psi\bU=(\tilde{\psi}_1, \ldots, \tilde{\psi}_N) \in V^N$ consists of $N$ pairwise
orthonormal  eigenfunctions of $A$, i.e. $ A \psi_k =\lambda_{k}
\psi_k $ for $k=1, \ldots, N$; in this case, there holds $\cJ(\Psi) = \sum_{k=1}^N \lambda_{k} $.
The minimum of $\cJ$ is attained if and only if the corresponding eigenvalues  $\lambda_k $, $ {k=1 , \ldots , N }$ are the $N$
lowest eigenvalues. This minimum is unique up to orthogonal
transformations if there is a gap ${ \lambda_{N+1} - \lambda_{N} >0
}$, so that in this case, the minimizers $\Psi= \mbox { argmin } \cJ$ are exactly the bases of the unique invariant
subspace spanned by the eigenvectors according to the $N$ lowest
eigenvalues.
\end{proposition}
%
%

\subsection{Comparison of direct minimization and\\ self consistent
iteration} 
Self consistent iteration consists of fixing the Fock operator $F^{(n)} = F_{\Phi^{(n)}}$ for each iterate $\Phi^{(n)}$; the simplified problem is then  solved in an inner iteration loo for $A = F^{(n)}$; the solution $\Phi$ defines the next iterate $\Phi^{(n+1)}$ of the outer iteration, by which the Fock operator is then updated to form $F^{(n+1)}$, defining the simplified problem for the next iteration step. For the solution of the inner problems with a fixed Fock operator, Proposition \ref{prop:localminima1} from the last section applies and the algorithms presented in the next section
can be used. Self consistent iteration is faced with
 convergence problems though, which can be remedied by advanced techniques: 
With an appropriate choice of the
update, the ODA-optimal damping algorithm \cite{lebris2}, convergence
can be guaranteed. \\
Direct minimization corresponds to the treatment of the nonlinear problem \eref{nonlinminiprob} for the Hartree-Fock or Kohn-Sham functional with the gradient algorithm 1 from the next section. Direct minimization thus differs from the self consistent iteration only in that the Fock operator is updated after each inner iteration step.
 Therefore, direct minimization is
preferable if 
the update of the Fock operator is sufficiently cheap. This is mostly the case for Gaussians
but not for the plane wave or wavelet basis or finite differences.

\section{Algorithms for minimization}\label{algorithsec}

In this section we will introduce three related algorithms to tackle the minimization problem \eref{nonlinminiprob} in a rather general form. Their convergence properties will be analysed in the next section.

\subsection{Gradient and projected gradient algorithm}

We will consider a gradient algorithm for the constrained
minimization problem; the motivation for this is given by the following related formulation (cf. \cite{lubich} for this concept):
With an initial guess $ \Phi^{(0)} \in \cV $,
$ [\Phi^{(0)}] \in \cG$,
 the gradient flow on $ \cV $, resp. $\cG$ is given by the differential
\begin{equation}\label{eq:ADE}
    \langle\langle \frac{d \Phi(t)}{dt}- \cJ'(\Phi(t)) , ~\delta \Phi   \rangle\rangle
= 0 ~~~ \forall \delta \Phi  \in \mathcal{T}_{[\Phi(t)]}\cG.\end{equation} 

Using the fact that $ \mathcal{I} -
\mathcal{D}_{[\Phi ]} $ is projecting onto the tangent space $
\mathcal{T}_{[\Phi]}\cG $, this algebraic differential initial value
problem can be rewritten by an ordinary initial value problem for
the gradient flow on $ \cV$,
\begin{equation}\label{eq:projcted-flow}
\frac{d}{dt} \Phi (t)   ~  =~ ( \mathcal{I} - \mathcal{D}_{[\Phi(t)]} )
\cJ'([\Phi(t) ]) ,~~~ \Phi (0) = \Phi^{(0)},
\end{equation}

or, equivalently,
\begin{equation}\label{eq:projcted-flow2}
\frac{d}{dt} \Phi (t)     = [ \cJ'  ,  \mathcal{D}_{[\Phi(t)]} ]
 (\Phi (t)) ,~~~ \Phi (0) = \Phi^{(0)},
\end{equation}
where the bracket $[.,..]$ denotes the usual commutator. Denoting by $(\cJ'(\Phi(t)))_i$ the $i$-th component of the gradient $\cJ'(\Phi(t))$ and letting
$\Lambda(t) := \big(\skapro{ (\cJ'(\Phi(t)))_i, \varphi_j(t) }\big)_{i,j=1}^N$, we obtain
the identification
\beqnr \cJ'(\Phi(t)) - \Phi(t)\Lambda(t) 
~ =~ [\cJ' , \mathcal{D}_{[\Phi(t) ]} ](\Phi(t)), \label{commutatorrel}\eeqnr
 which we will make use of later.

There holds
$\frac{d \Phi(t)}{dt} \to 0$ for $t \to \infty$, so we are looking for the fixed point of this flow $ \Psi = \lim_{t\to
\infty } \Phi (t)$
rather than its trajectory. Equation \eref{eq:projcted-flow2} suggests the projected gradient type
algorithms presented below. In algorithm 1, corresponding to an Euler procedure for the differential equation \eref{eq:ADE}, 
the gradient at a certain point $\Phi(t)$ is kept fixed (and being preconditioned) for non-differential stepsize, so that the manifold is left in each iteration step. Therefore, a projection on the admitted set is performed in each iteration step.\\
Note also that the role of the preconditioners
 $ \cB_{n}^{-1}$ is crucial, see the remarks following algorithm 1.\smallskip

\textit{\textbf{Algorithm 1: Projected Gradient Descent}}\vspace*{0.5mm} 
\footnotesize
\hrule \smallskip\noindent
\textit{\textbf{Require:} Initial iterate $\Phi^{(0)} \in V$;\\
\phantom{Require: }~evaluation of $\cJ'(\Phi^{(n)})$ and of preconditioner(s) $B_n^{-1}$}  (see comments below) \\
\textit{\textbf{Iteration:} \\
for $n= 0, 1, \ldots$ do\\
$\phantom{o}\quad (1)$ Update $\Lambda^{(n)} := \skapro{\cJ'(\Phi^{(n)}), \Phi^{(n)} } \in \R^{N \times N}, $\\
 $\phantom{o}\quad (2)$ Let $\hat{\Phi}^{(n+1)} := \Phi^{(n)} -  \cB^{-1}_n(\cJ'(\Phi^{(n)}) -\Phi^{(n)}\Lambda^{(n)}),$ \\
$\phantom{o}\quad$ \phantom{(2) }$\big(=\Phi^{(n)} -  \cB^{-1}_n(\cA_{\Phi^{(n)}}\Phi^{(n)} -\Phi^{(n)}\Lambda^{(n)})$ for the case that $\cJ'(\Phi) = \cA_{\Phi}\Phi.\big)$\\
 $\phantom{o}\quad (3)$ Let $ \Phi^{(n+1)} = P\hat{\Phi}^{(n+1)}$ by projection $P$ onto
$ \mathcal{V}$ resp. $\mathcal G$\\
endfor}\smallskip
\hrule
\normalsize
\bigskip\bigskip

Some remarks about this algorithm are in order. First of all, note that if Algorithm 1 is applied to the ansatzes in electronic structure calculation as portrayed in section \ref{schrosec}, the gradient $\cJ' (\Phi) $ is given by $\cJ' (\Phi) = \cA_{\Phi}\Phi$ with $A_{\Phi}$ the Fock- or Kohn-Sham operator or a fixed operator $A_{\Phi} = A$  for the simplified problem. Therefore,
$(\cJ'(\Phi^{(n)} -\Phi^{(n)}\Lambda^{(n)})_i = A_{\Phi^{(n)}}\phi_i^{(n)} - \sum_{j=1}^N  \ll A_{\Phi^{(n)}}\phi_i^{(n)}, \phi_j^{(n)}\rr\phi_j^{(n)} $ is the usual ``subspace residual'' of the iterate $\Phi^{(n)}$, which is a crucial fact for capping the complexity of the algorithm in section \ref{commentsec}.\\

Next, let us specify the role of the preconditioner $ \cB_n^{-1}$ used in each step. This preconditioner is induced (according to \eref{calliVN}) 
by an elliptic symmetric operator $B_n: V \to V'$, which
we require to be equivalent to
the norm on $H^1$  in the sense that
\begin{equation}\label{eq:eq-precond}
\Lska{B_n \varphi, \varphi} ~~ \sim~~
\norm{\varphi}^2_{H^1} ~~~\forall \varphi \in V = H^1(\R^3)\ .
\end{equation}

For example, one can use approximations of  the shifted Laplacian,
$B \approx \alpha (-\frac{1}{2} \Delta + C) $, as is done in the BigDFT project. This is also 
a suitable choice when dealing with plane wave ansatz functions using
advantages of FFT, or a multi-level preconditioner
 if one has finite differences,
finite elements or multi-scale functions like wavelets \cite{hackbusch, becker, goedecker, ariaswave}.

For the simplified problem, the choice $ B^{-1} = \alpha A^{-1} $
corresponds to a variant of simultaneous inverse iteration. The
choice  
$$ B|_{ V_0^{\perp} := \{v | \langle v , \varphi_i^{(n)} \rangle = 0 
\, \forall i=1, \ldots , N \}}
 = \alpha ( A    - \lambda_j^{(n)} I  )|_{ V_0^{\perp} := \{v | \langle v ,
\varphi_i^{(n)} \rangle = 0 \, \forall i=1, \ldots , N \}} $$
 corresponds
to a simultaneous Jacobi-Davidson iteration.

To guarantee convergence of the algorithm, the preconditioner $B$ chosen according to the guidelines above also has to be properly scaled by a factor $\alpha > 0$, cf. Lemma \ref{lem1}. The optimal choice of $\alpha$ is provided by minimizing the corresponding functional over
$\mbox{span } \{ \Phi^{(n)} , \widehat{\Phi}^{(n+1)} \} $ (a line search over this space), which can be
done for the simplified problem without much additional effort. For
the Kohn-Sham energy functional, it will become prohibitively
expensive. However, line
search and subspace acceleration like DIIS \cite{diis}
will improve the convergence speed. Note that in this context, one might as well use 
different step sizes for every entry,
i.e. $\cB \Phi = (\alpha_1 B \varphi_1, \ldots, \alpha_N B \varphi_N)$.

Next, let us make a remark concerning the projection onto $ \cG$. It
only has to satisfy $ \mbox{span } \{\varphi_i^{(n+1)} :   1 \leq i \leq N \} = \mbox{span } \{
\widehat{\varphi}_i^{(n+1)} :   1 \leq i \leq N \} $. For this purpose
any orthogonalization of $ \{ \widehat{\varphi}_i^{(n+1)} :   1 \leq i
\leq N \} $ is admissible. For example, three favorable possibilities
which up to unitary transformations yield the same result are 
\begin{itemize}
  \item Gram-Schmidt orthogonalization,
    \item Diagonalization of the Gram matrix
    $ \mathbf{G} = ( \langle  \hat{\varphi}_i^{(n+1)} , \hat{\varphi}_j^{(n+1)}\rangle)_{i,j=1}^N
     $ by Cholesky factorization,
    \item (For the problems of section \ref{schrosec}, i.e. where $\cJ'(\Phi) = \cA_{\Phi}\Phi$:) \\ Diagonalisation of the matrix
    $ \mathbf{A}_{\Phi^{(n+1)}} := ( \langle A_{\Phi^{(n)}} \hat{\varphi}_i^{(n+1)} , \hat{\varphi}_j^{(n+1)} \rangle)_{i,j=1}^N
    $  by solving an $ N\times N$ eigenvalue problem.\unred
\end{itemize}

Parallel to the above algorithm, we consider the following variant in which the descent direction
is projected onto the tangent space $\cT_{[\Phi^{(n)}]}\cG$ in every iteration step.
It will play an important theorectical role considering convergence of the local exponential
parametrization, i.e. algorithm 3.

\textit{\textbf{Algorithm 2: Modified Projected Gradient Descent}}\vspace*{0.5mm} 
\footnotesize
\hrule \medskip\noindent
\textit{\textbf{Require:} see Algorithm 1} \\
\textit{\textbf{Iteration:} \\
for $n= 0, 1, \ldots$ do\\
$\phantom{o}\quad (1)$ Update $\Lambda^{(n)} := \skapro{\cJ'(\Phi^{(n)}), \Phi^{(n)} } \in \R^{N \times N}, $\\
$\phantom{o}\quad (2)$ Let $\hat{\Phi}^{(n+1)} := \Phi^{(n)} -  (\cI-\cD_{\Phi^{(n)}}) \cB^{-1}_n \big(\cJ'(\Phi^{(n)}) -\Phi^{(n)}\Lambda^{(n)}\big),$ \\
$\phantom{o}\quad$ \phantom{(2) }$\big(=\Phi^{(n)} -  \cB^{-1}_n(\cA_{\Phi^{(n)}}\Phi^{(n)} -\Phi^{(n)}\Lambda^{(n)})$ for the case that $\cJ'(\Phi) = \cA_{\Phi}\Phi.\big)$\\
$\phantom{o}\quad (3)$ Let $\Phi^{(n+1)} = P \hat{\Phi}^{(n+1)} $ by projection $P$ onto
$ \mathcal{V}$ resp. $\mathcal G$,\\
endfor} \medskip
\hrule
\normalsize
\medskip

Note again that the algorithms are given in a general form, where the preconditioner (or the corresponding parameter $\alpha$, e.g. obtained by a kind of line search) may be chosen in each iteration step. In our analysis,  we will consider a fixed preconditioner $\cB_n = \cB$ in every iteration step, for which we will show linear convergence without further line search invoked. Thus, our analysis is in a way a worst case analysis for the algorithms under consideration. See also section \ref{commentsec} for improvements on the speed of convergence.

\subsection{Exponential parametrization}

Instead of projecting the iterate $\Phi^{(n)}$ onto the Grassmann
manifold $\cG$
in every iteration step, we will now develop an algorithm in which
the iterates remain on the manifold without further projection.
This will be achieved by following geodesic paths on the manifold
instead of straight
lines in Euclidean space, which has the advantage that during our
calculations we
do not leave the constraining set at any time so that no orthonormalization process is required.
To apply the result of proposition \ref{geodesics}, we will for this algorithm limit our treatment to the case where $\cJ'(\Phi) = \cA_{\Phi}\Phi$ is given by a linear operator (see the discussion after algorithm 1).\\
Recall that a
geodesic is a curve $c$ on a manifold with vanishing second covariant
derivative,
i.e.
\beqnr \frac{\nabla}{dt} \dot{c}(t) ~:= ~\pi_{c(t)} \ddot{c}(t) ~=~ 0 ~~~\mbox{ for
all } t,\label{covariant}\eeqnr
where $\pi_{c(t)}$ denotes the projection onto the tangent space at the
point $c(t)$.

\begin{proposition} \label{geodesics}
For any operator $X: V \to V$ for which $\XX \Phi \in \mathcal{T}_{[\Phi]} \mathcal{G}$ (where as always, $\XX$ is defined by $X$ by \eref{calliVN}),
the antisymmetric operator
\beqnr \hat{X} = (I-D_{\Phi})X D_{\Phi} - D_{\Phi} X^{\dagger} (I-D_{\Phi}),\label{expoantisym}\eeqnr
satifies  $\hat{\XX} \Phi = \XX \Phi$, and $c(t) := \exp(t\hat{\XX}) \Phi$ is a geodesic in 
$\mathcal{G}$
emanating from point $\Phi$
with direction $\dot{c}(0)= \XX \Phi$.
\end{proposition}
\begin{proof}
The proof is straightforward; application of the projection equation yields
$
\left( \frac{\nabla}{dt} \dot{c}(t) \right)
 = (\cI - \cD_{c(t)}) \ddot{c}(t) = \mathbf{0}.$
\end{proof}

If we now let, for any iterate $\Phi^{(n)}$, \beqnr X^{(n)} = (I-D_{\Phi^{(n)}})B^{-1}(I-D_{\Phi^{(n)}})A_{\Phi^{(n)}}, \label{expoop}\eeqnr the curve \beq c(t)&:=& \exp (- t \hat{\XX}^{(n)}) \Phi^{(n)}\eeq with $\hat{\XX}^{(n)}$ from \eref{expoantisym}
is by the previous Lemma a geodesic in $\cG$ with direction
\begin{eqnarray*}
\dot{c}(0)~=~- (I-D_{\Phi^{(n)}})B^{-1}(I-D_{\Phi^{(n)}})A_{\Phi^{(n)}} \Phi^{(n)}
\end{eqnarray*}
which equals the (preconditioned) descent direction of the projected gradient descent
algorithm of the preceding section.
If we now choose the next iterate as a point on this geodesic, we get the
following algorithm:\medskip

\textit{\textbf{Algorithm 3: Preconditioned exponential parametrization}}\vspace*{0.5mm}
\footnotesize
\hrule \medskip\noindent
\textit{\textbf{Require:}  see Algorithm 1} \\
\textit{\textbf{Iteration:}}\\
\textit{for} $n= 0, 1, \ldots$ \textit{do}\\
$\phantom{oo} \rhd$ \textit{Follow a geodesic path on the Grassmann manifold with stepsize}
$\alpha$, 
$$ \Phi^{(n+1)} := \exp (- \alpha \hat{\XX}^{(n)} )\Phi^{(n)} \makebox{(with $\XX$ from \eref{expoop} and $\hat{\XX}$ defined by  $\XX$ via \eref{expoantisym})}$$
\textit{endfor} \medskip
\hrule
\normalsize
\bigskip

Note that a similar algorithm,
Conjugate Gradient on the Grassman Manifold, has already
been introduced in \cite{arias}, page 327. That paper also included numerical
tests for a model system. The algorithm was also tested for electronic structure applications very different from those of the BigDFT program in \cite{raczkow}. A similar approach using the density matrix representation for electronic structure problems  was also proposed in \cite{shao}, where the authors move along the geodesics in a gradient resp. Newton method direction without preconditioning. \\
Like in this work, the stepsize $\alpha$ may be calculated in each iteration step using
line search algorithms like backtracking linesearch or quadratic approximations to
the energy term \cite{nocedal}. These often time consuming line searches may be omitted though if we choose a
suitable preconditioner
$B = B_n$ and set the stepsize $\alpha = 1$ once and for all.\\
The efficiency of this algorithm strongly depends on the computation
of matrix exponentials needed to follow geodesic paths on the Grassmann
manifold.
A variety of methods can be found in \cite{dubiousways}, see also \cite{saad} for an analysis of selected methods.
For some of these algorithms, there exist powerful implementations like
the software package Expokit
\cite{expokit}, which contain both Matlab and Fortran code thus supplying a convenient 
tool for numerical experiments.

\section{Convergence results}\label{convsec}

\subsection{Assumptions, error measures and main result}

In this section, we will show linear convergence of the algorithms of the last section under the ellipticity assumption \ref{nonlinasselli} given below.
Additional results we give include the equivalence of the error of $\Phi$, measured in a norm on $V$, and the error of the gradient residual $(\cI-\cD)\cJ'(\Phi)$, and quadratic reduction of the energy error $\cJ(\Phi^{(n)})-\cJ(\Psi)$.\\
Recall that in our framework introduced in section \ref{theoframesec}, we kept  the freedom of choice to either use $V := H^1 = H^1(\R^3)$, equipped with an inner product equivalent to the $H^1$
inner product $\ll .,.. \rr_{H^1}$, for analysing the original equations,
or to use $V = V_h \subset H^1 $ as a finite dimensional subspace for a
corresponding Galerkin discretization of these equations. In practice, our iteration scheme
is only applied to the discretized equations. However, the convergence
estimates obtained will be uniform with respect to the discretization
parameters.
The main ingredient our analysis is based on is the following condition imposed on the functional $\cJ$, cf. section \ref{ministiefel}:
\begin{assumption}\label{nonlinasselli}
Let $\Psi$ a minimizer of \eref{nonlinminiprob}. The Hessian $ \cL^{(2,\Psi)} (
\Psi, \Lambda) : V^N \to (V')^N $ of the Lagrangian $\cL(
\Psi, \Lambda)$ (given by \eref{lagrangian}), where the derivatives are taken with respect to $\Psi$,
is assumed to be $V^N$-elliptic on the tangent space, i.e. there is $\gamma > 0$ so that
\beqnr
\langle\langle \cL^{(2,\Psi)} ( \Psi , \Lambda) \delta \Phi ~,~ \delta \Phi   \rangle\rangle~~\geq~~
\gamma ~\| \delta \Phi \|_{V^N}^2, ~~~~ \makebox{for all}~~~ \delta \Phi \in \mathcal{T}_{[\Psi]}\cG.
\label{lagrangeelli}
\eeqnr
\end{assumption}
Note again that for Hartree-Fock calculations, verification of $\cL^{(2,\Psi)}( \Psi , \Lambda) > 0$ on $\cT_{[\Psi]}\cG$ already implies $\cL^{(2,\Psi)}( \Psi , \Lambda)$, cf \cite{maday}.

From section \ref{grasstiefel}, we recall that $ \cL^{(2,\Psi)} (
\Psi, \Lambda)\Phi = \cJ''(\Psi)\Phi - \Phi\Lambda $, so that \eref{lagrangeelli} is verified if and only if
\beqnr
\langle\langle \cJ'' (\Psi) \delta \Phi - \delta \Phi \Lambda ~,~ \delta \Phi   \rangle\rangle~~\geq~~
\gamma ~\| \delta \Phi \|_{V^N}^2, ~~~~ \makebox{for all}~~~ \delta \Phi \in \mathcal{T}_{[\Psi]}\cG \label{expliccond}
\eeqnr holds, where $\Lambda = (\ll (\cJ'(\Psi))_j,\psi_i\rr)_{i,j=1}^N$ as above.
From the present state of Hartree-Fock theory, it is not possible to decide whether this condition is true in general; the same applies to DFT theory. For the simplified problem, the condition holds if the operator $A$ fulfils the conditions of the following lemma.

\begin{Lemma}\label{linearprobelli}
Let $A: V \to V', \psi
\mapsto A \psi$ a bounded symmetric
operator, such that $A$ has $N$ lowest eigenvalues $\lambda_1 \leq \ldots \leq \lambda_N$
satisfying the gap condition  \beqnr \lambda_N ~~<~~ \inf \{\lambda ~|~ \lambda \in
 \sigma (A) \backslash \{  \lambda_1, \ldots , \lambda_N \}\}.
\label{gapcondition} \eeqnr 
Then assumption \ref{nonlinasselli} holds for the simplified problem \eref{miniprob}.
\end{Lemma}

\begin{proof} We estimate the two terms of \eref{expliccond} separately. Let us denote $\lambda = \inf \{\lambda ~|~ \lambda \in
 \sigma (A) \backslash \{  \lambda_1, \ldots , \lambda_N \}\}.$ To the first term, the Courant-Fisher theorem (\cite{minmax}) applies componentwise to give the estimate $ \llll \mathcal{A}
\delta\Phi , \delta\Phi \rrrr \geq \lambda || \delta\Phi||_{V^N}^2$.
For the second, choosing $\bU = (u_{i,j})_{i,j=1} ^N \in O(N)$ so that $\bU^T \Lambda \bU = \diag(\lambda_i)_{i=1} ^N,$ where $\lambda_i$ are the lowest $N$ eigenvalues of $A$, gives 
\beq 
\llll  \delta\Phi \Lambda  , \delta\Phi\rrr~~&=&~~ \llll \delta\Phi (\bU \bU^T \Lambda \bU \bU^T) , \delta\Phi\rrr 
~~~:=~~ \sum_{i=1}^N \ll \sum_{j=1}^N u_{j,i} \lambda_{j} \delta\varphi_j, \sum_{k=1}^Nu_{k,i}\delta\varphi_k \rr \\&=&~~   \sum_{j, k =1}^N  \lambda_{j} \delta_{j,k} \ll \delta\varphi_j, \delta\varphi_k \rr ~~~~\leq~~ \lambda_N || \delta\Phi||_{V^N}^2.\\
\eeq 
 so that $\cL^{(2)}(\Psi, \Lambda) $ is elliptic on  $\cT_{[\Psi ]} \mathcal{G}$ by the gap condition \eref{gapcondition}.
\end{proof}

To formulate our main convergence result, we now introduce a norm $||.||_{V^N}$ on the space $V^N$, which will be equivalent to the $(H^1)^N$-norm but more convenient for our proof of convergence. We will then state our convergence result in terms of these error measures.

\begin{Lemma}\label{lem:def-b}
Let $B: V \to V'$ the preconditioning mapping introduced in section \ref{algorithsec}, so that in particular, $B$ is symmetric and the spectral equivalence
$$ \vartheta||x||_{H^1}^2 ~~\leq ~~ \ll B x, x \rr  ~~\leq ~~ \Theta ||x||_{H^1}^2$$
holds for some $ 0 < \vartheta \leq \Theta$ and all $x \in V$.
Let us consider the mapping
\begin{equation}\label{eq:def-hatb}
\hat{B}^{-1} : V' \to V \ , ~~
    \hat{B}^{-1} := ( I - D ) B^{-1} (I - D ) + D ,
\end{equation}
where $D = D_{\Psi}$ projects onto the sought subspace. Then the inverse $ \hat{B} $ satisfies $
 \langle \hat{B}  \varphi , \psi  \rangle  =  \langle \varphi , \hat{B}  \psi \rangle
$ for all $\varphi, \psi \in V$, and for the induced $\hat{B}$-norm $||.||_{\hat{B}}$ on $V$ there holds
$$
 \langle \hat{B}  \varphi , \varphi  \rangle ~~ \sim~~ \|\varphi \|^2_{H^1} .
$$\end{Lemma}

Using the notation \eref{calliVN}, a norm  on $V^N$ is now induced by the $||.||_{\hat{B}}$-norm by
\begin{equation}\label{eg:defB}
    \| \Phi \|_{V^N}^2  :=  \langle \langle  \widehat{\mathcal{B}} \Phi , \Phi
    \rangle\rangle  \ .
\end{equation}
Note that this norm, as any norm defined on $V^N$ in the above fashion, is invariant under the orthogonal group of $\R^{N \times N}$ in the sense that \beqnr \| \Phi\mathbf{U} \|_{V^N} = \| \Phi \|_{V^N} \label{invarianceofnorm}\eeqnr
for all $\bU \in O(N)$.
In the Grassmann manifold, we measure the error between
$[\Phi_{(1)}], [\Phi_{(2)}] \in \mathcal G$ by a related metric $d$ given by
$$d(~{[\Phi_{(1)}]} , [\Phi_{(2)}]~) ~:=~
\inf_{\mathbf{U} \in O(N)} \| \Phi_{(1)} -\Phi_{(2)} \mathbf{U}
\|_{ V^N}.
$$
If $ [\Phi_{(2)}] $ is sufficiently close to $
[\Phi_{(1)} ] \in \mathcal{G}$ it follows from Lemma \ref{lem:projdiffer} that this measure given by $d$
is equivalent to the expression
\begin{equation}\label{eq:diff-norm}
 \|  ( \cI - \mathcal{D}_{\Phi_{(1)}}) \Phi_{(2)} \|_{V^N},
\end{equation}
in which we used the $L_2$-orthogonal projector $\mathcal{D}_{\Phi_{(1)}}$ onto the subspace spanned by $\Phi_{(1)}$.
In the following, let us use the abbreviation $D = D_{\Psi}$ for the projector on the sought subspace, whereever no confusion can arise. An equivalent error measure for the deviation of $\Phi \in \cV$ from the sought element $\Psi \in \cV$ is then given by the expression
\begin{equation} \label{eq:b-error}\|  ( \cI - \mathcal{D}) \Phi\|_{V^N},\end{equation}
which will be used in the sequel. In terms of this notation, our main convergence result is the following.

\begin{theorem}\label{mainconv}
Under the ellipticity assumption \eref{lagrangeelli}, the following holds for any of the three algorithms formulated in section \ref{algorithsec}: For $\Phi^{(0)} \in U_{\delta}(\Psi)$ sufficiently close to $\Psi$, there is a constant
$\chi < 1$ such that for all $n \in \N_0,$
\beqnr \norm{(\cI - \cD) \Phi^{(n+1)}}_{V^N} ~~~\leq~~~
\chi \cdot \norm{(\cI - \cD) \Phi^{(n)}
}_{V^N}.\label{mainconvexpr} \eeqnr
\end{theorem}


The rest of this section will be mainly dedicated to the proof of this theorem. 
For the sake of clarity, let us first sketch the proof to be performed: We will exploit the fact that the iteration mapping can be written in the form  $\Phi^{(n)} \mapsto \Phi^{(n)} - \cB^{-1} (\cI - \cD_{\Phi^{(n)}}) \cJ'( \Phi^{(n)})$ and is thus a perturbation of the mapping $ \Phi^{(n)} \mapsto \Phi^{(n)} - \cB^{-1}  (\cI - \cD_{\Psi}) \cJ'(\Phi^{(n)})$.
 The estimate then splits in two main parts: The first will be a linear part incorporating the Hessian of the Lagrangian and the task will be to show that application of this linear part to an iterate $\Phi^{(n)} \in \cG$ indeed reduces its error in the tangent space of $\Psi$ (as defined by \eref{eq:b-error}); here, our ellipticity assumption enters as main ingredient. The second part consists of showing that the remaining perturbation terms (including those resulting from projection on the manifold) are of higher order and thus asymptotically neglectable; the main lemmas entering are Lemma \ref{lem:projdiffer} above and Lemma \ref{lem:hat-p} to be proven below. 

\subsection{Ellipticity on the tangent space}

In this section, we will first formulate a rather general result about how ellipticity on subspaces can be used to construct a contraction on these spaces and then specialize this to the tangent space at the solution $\Psi$ and assumption \ref{nonlinasselli} in the subsequent corollary. Finally, we will then prove that our assumption \ref{lagrangeelli} entering here is indeed true for the simplified problem \eref{miniprob}.

\begin{Lemma}\label{lem1}
Let $ W \subset G \subset W' $ a Gelfand triple, $ U \subset W$ a
closed subspace of $W$
 and    $S, T': W
\to W'$ two bounded elliptic operators, symmetric with respect
to the $G$-inner product $\ll.,..\rr_G$, satisfying
 \beqnr&&
\gamma ||x||_W^2 ~~\leq ~~ \ll S x, x \rr_G  ~~\leq ~~ \Gamma ||x||_W^2,\\
 \makebox {and}~~~~&& \vartheta||x||_W^2 ~~\leq ~~ \ll T' x, x \rr_G  ~~\leq ~~ \Theta ||x||_W^2 \label{abschAB}
\eeqnr for all $x \in U$. Moreover, let $S, T'$ both map the subspace $U$ to itself. Then there exists a scaled variant $T =
\alpha T'$, where $\alpha > 0$, and a constant $\beta < 1$ for which 
\beqnr || (I - T^{-1}S ) x||_T ~~\leq ~~ \beta~ ||x||_T.
\label{itmatcontr},\eeqnr 
 for all $ x \in U$, where $ ||x||_T^2 :=  \ll T x, x \rr_G $ is the inner product induced by $T$. 
\end{Lemma}

\begin{proof}
It is easy to verify that for $\beta := (\Gamma \Theta - \gamma \vartheta)/(\Gamma \Theta +
\gamma \vartheta) < 1$
and
$\alpha := \frac 12 (\Gamma/\vartheta + \gamma / \Theta)$
there holds\beqnr |\ll (I - T^{-1}S ) x, x \rr_T| ~~\leq ~~ \beta~ ||x||_T^2~~~\mbox{ for all } x \in U.\eeqnr 
Due to the symmetry of $T, S$ as mappings $U \to U$, the result \eref{itmatcontr} follows. \end{proof}

Let $\lambda_i, i =  1, \ldots, N$ the lowest eigenvalues of $A$, $\psi_i, i = 1, \ldots, N$, the corresponding eigenfunctions, and \beqnr V_0 = \mbox{span } \{ \psi_i : i = 1,
\ldots , N \} \label{v0defi}\eeqnr

By Lemma \ref{stiefeltangentspace}, there holds $(V_0^{\bot})^N = \mathcal{T}_{[\Psi ]} \mathcal{G} ,$ where $\Psi = (\psi_1, \ldots, \psi_N).$
The following corollary is the main result needed for estimation of the linear part of the iteration scheme.

\begin{corollary}\label{iterationconv}

Let $\cJ$ fulfil the ellipticity condition \eref{lagrangeelli} and $B': V \to V'$ a symmetric operator that fulfils \eref{abschAB} with $T' = B'$. Then there exists a 
scaled variant $B = \alpha B'$, where $\alpha > 0$, for which for any $ \delta\Phi \in \mathcal{T}_{[\Psi ]} \mathcal{G} $
there holds
$$ \|  \delta\Phi - \hat{\mathcal{B}}^{-1}(\cI - \cD)\cL^{(2, \Psi)}(\Psi,  \Lambda) \delta\Phi   \|_{V^N} ~~ \leq~~ \beta~ \| \delta\Phi
\|_{V^N},
$$
where $\beta <1 $ and $\hat{B}$ is defined by $B$ via \eref{eq:def-hatb}.
\end{corollary}

\begin{proof} Note that the restriction of $\hat{B'}$ is a symmetric operator $V_0^{\perp} \to V_0^{\perp} $, so that the same holds for the extension $\hat{\cB'}$ as mapping $\cT_{[\Psi ]} \mathcal{G} \to \cT_{[\Psi ]} \mathcal{G}.$ 
$(\cI - \cD)\cL^{(2, \Psi)}$ also maps $V_0^{\perp} \to V_0^{\perp} $ symmetricly, so Lemma \ref{lem1} applies.
\end{proof}

\subsection{Residuals and projection on the manifold}

For the subsequent analysis, the following result will be useful. It also shows that the ``residual'' $(\cI-\cD_{\Phi^{(n)}})\cJ'(\Phi^{(n)})$
may be utilized for practical purposes to estimate the norm of the error $(I-D)\Phi^{(n)}$.

\begin{Lemma}
For $\delta$ sufficiently small and $||(\cI-\cD)\Phi^{(n)}||_{\hat{B}} < \delta$, there are constants $c, C>0$ such that
\beqnr
c||(\cI-\cD)\Phi^{(n)}||_{V^N} ~~\leq~~ ||(\cI-\cD_{\Phi^{(n)}})\cJ'(\Phi^{(n)})||_{(V^N)'} ~~\leq~~ C||(\cI-\cD)\Phi^{(n)}||_{V^N} .\label{resifirst}
\eeqnr
An analougeous result holds for gradient error $||(\cI-\cD)\cJ'(\Phi^{(n)})||_{(V^N)'}$.
\end{Lemma}
\begin{proof}
Let us choose $\bar{\Psi} \in [\Psi]$ according to Lemma \ref{lem:projdiffer} (applied to $\Phi = \Phi^{(n)}$). 
Letting $\Delta \Psi :=\Phi^{(n)} - \bar{\Psi} ,$ there holds by linearization and Lemma \ref{lem:projdiffer} (recall that we let $D=D_{\Psi}$)
\beq(\cI-\cD_{\Phi^{(n)}})\cJ'(\Phi^{(n)}) &=& (\cI-\cD)\cJ'(\Psi)~ +~ (\cI-\cD)\cL^{(2,\Psi)}(\bar{\Psi},\Lambda) \Delta\bar{\Psi} ~+~ \cO(||(\cI-\cD)\Phi^{(n)}||^2_{V^N}\\
&=& (\cI-\cD)\cL^{(2,\Psi)}(\bar{\Psi},\Lambda)(\cI-\cD)\Phi^{(n)} + \cO(||(\cI-\cD)\Phi^{(n)}||^2_{V^N})\eeq
By assumption \ref{nonlinasselli}, $||(\cI-\cD)\cL^{(2,\Psi)}(\Psi,\Lambda)(\cI-\cD)\Phi^{(n)}||_{(V^N)'} \sim ||(\cI-\cD)\Phi^{(n)}||_{V^N}$, from which the assertion follows. The assertion for  $||(\cI-\cD)\cJ'(\Phi^{(n)})||_{(V^N)'}$ follows from the same reasoning by replacing $\cL^{(2,\Psi)}(\Psi,\Lambda)$ by $\cJ''(\Psi)$ in the above.
\end{proof}

The last ingredient for our proof of convergence is following lemma which will imply that the projection following each application of the iteration mapping does not destroy the asymptotic linear convergence.

\begin{Lemma}\label{lem:hat-p}
Let $\hat{\Phi}^{(n+1)} = (\hat{\phi}_1, \ldots , \hat{\phi}_N)$ the intermediate iterates as resulting from iteration step (2) in algorithm 1 or 2, respectively.
For any orthonormal set $\Phi \in \cV$ fulfilling $\spn[\Phi] = \spn[\hat{\Phi}^{(n+1)}]$, its error 
deviates from that of $\hat{\Phi}^{(n+1)}$ only by quadratic error term:
\beqnr
||(\cI-\cD)\Phi ||_{V^N} ~~=~~ ||(\cI-\cD)\hat{\Phi}^{(n+1)}||_{V^N}~+~\cO(||(\cI-\cD)\hat{\Phi}^{(n)}||_{V^N}^2)   \label{projesti}
\eeqnr 
\end{Lemma}

\begin{proof} 
First of all, note that if \eref{projesti} holds for one orthonormal set $\Phi$ with $\spn[\Phi] = \spn[\hat{\Phi}^{(n+1)}]$, it holds for any other orthonormal set $\tilde{\Phi}$  with $\spn[\tilde{\Phi}] = \spn[\hat{\Phi}^{(n+1)}]$ because $||(\cI-\cD)\Phi\bU||_{V^N} = ||(\cI-\cD)\Phi||_{V^N}$ for all orthonormal $\bU \in O(N)$. Therefore, we will show \eref{projesti} for $\Phi = (\varphi_1, \ldots, \varphi_N)$ yielded from $\hat{\Phi}^{(n+1)}$ by the Gram-Schmidt orthonormalization procedure. 
Denote $\hat{\varphi}_i = \varphi^{(n)}_i+ r_i^{(n)}$, where for $s_i^{(n)} = B^{-1}\big(\cI-\cD_{\Phi^{(n)}})\cJ'(\Phi^{(n)}),$ we set $r_i^{(n)} = s_i^{(n)}$ or $r_i^{(n)} = (I-D_{\Phi^{(n)}})s_i^{(n)}$ for algorithm 1 or 2, respectively. From the previous lemma, we get in particular that $||r_i^{(n)}||_{V} \lesssim ||(I-D)\phi_i^{(n)}||_{V}$ for both cases (remember that $D = D_{\Psi}$). 
With the Gram-Schmidt procedure given by $\varphi'_k = \hat{\varphi}_k - \sum_{j<i}\ll \hat{\varphi}_k,\varphi_j\rr\varphi_j$, $\varphi_k = \varphi'_k/||\varphi'_k||$,  the lemma is now proven by verifying that in each of the inner products involved, there occurs at least one residual $||r_i^{(n)}||$; and that, on top of this, for the correction directions $\varphi_j$ there holds $(I-D)\varphi_j' = \cO(||(\cI-\cD)\Phi^{(n)}||_{V^N}) + \cO(\sum_{i<k}||r_i^{(n)}||_{V^N})= \cO(||(\cI-\cD)\Phi^{(n)}||_{V^N})$. Therefore, the correction terms are of $\cO(||(\cI-\cD)\hat{\Phi}^{(n)}||_{V^N}^2)$, thus proving $\varphi'_k - \hat{\varphi}_k = \cO(||(I-D)\Phi||_{V^N}^2).$ It is easy to verify that the normalization of $\varphi'_k$ only adds another quadratic term, so the result follows.
\end{proof}
\subsection{Proof of Convergence}\label{convsubsec}

To prove \eref{mainconvexpr} for Algorithm 1, we define $\cF(\Phi) = \Phi - \mathcal{B}^{-1}
(\cI-  \mathcal{D}_{\Phi} )  \cJ'(\Phi)$, so that $\Phi^{(n+1)} = P (\cF(\Phi^{(n)})),$
where $P$ is a projection on the Grassmann manifold for which $ [P (\cF(\Phi^{(n)}))] = [\cF(\Phi^{(n)})]$.
For fixed $n$, let us choose $\bar{\Psi} \in \spn[\Psi]$ according to Lemma \ref{lem:projdiffer}, so that, using the abbreviation $\mathcal{D}:= \mathcal{D}_{\Psi }$,
\beqnr  \bar{\Psi} - \Phi^{(n)}~~&=&~~ (\cI - \cD)  \Phi^{(n)} + \cO(||(\cI - \cD)  \Phi^{(n)}||_{L^N_2}^2)\\
& \leq &~~ (\cI - \cD)  \Phi^{(n)} + \cO(||(\cI - \cD)  \Phi^{(n)}||_{V^N}^2) \label{projectcite}\eeqnr
Introducing $\Delta \Psi :=\Phi^{(n)} - \bar{\Psi} ,$ there follows by linearization
\beqnr &&\|( \cI- \mathcal{D}) \Phi^{(n+1)}   \|_{V^N} \\&\stackrel{\textsf{Lemma}~\ref{lem:hat-p}}{=}&
 \| ( \cI- \mathcal{D})  \cF(\Phi^{(n)})  \|_{V^N} + \mathcal{O} ( \|  ( \cI- \mathcal{D}) \Phi^{(n)} \|_{V^N}^2 )\\
&=&  \|  ( \cI- \mathcal{D})  \cF( \bar{\Psi})   +
 ( \cI- \mathcal{D}) \cF'(\bar{\Psi}) \Delta\Psi\|_{V^N} + \mathcal{O} ( \|  ( \cI- \mathcal{D})\Phi^{(n)} \|_{V^N}^2 )\\
&=&  \|  ( \cI- \mathcal{D})  \cF'(\bar{\Psi}) (\cI - \cD)\Phi^{(n)} \|_{V^N}
 + \mathcal{O} ( \|  ( \cI- \mathcal{D})\Phi^{(n)} \|_{V^N}^2 )\\
&=&  \|  ( \cI- \mathcal{D})\big (\cI - \cB^{-1}(\cI-\cD)\cL^{(2, \Psi)}(\bar{\Psi}, \Lambda)\big) (\cI - \cD)\Phi^{(n)} \|_{V^N} \label{estimline}\\ && + \mathcal{O} ( \|  ( \cI- \mathcal{D})\Phi^{(n)} \|_{V^N}^2 ) \notag \eeqnr
where we have used \eref{projectcite} and the fact that $( \cI- \mathcal{D})  \cF( \bar{\Psi})$ is zero.
The proof is now finished by noticing that
\beq 
&& (\cI - \cD)  \Big(\cI -  \cB^{-1}(\cI - \cD) \cL^{(2, \Psi)}(\bar{\Psi}, \Lambda)\Big) (\cI - \cD)\Psi\\
&=& \Big(\cI - \hat{\cB}^{-1}(\cI - \cD)\cL^{(2, \Psi)}(\bar{\Psi}, \Lambda)\Big) (\cI - \cD)\Psi,
\eeq 
so that corollary \ref{iterationconv} applies to give 
$$\|( \cI- \mathcal{D}) \Phi^{(n+1)}   \|_{V^N} ~\leq ~ \vartheta  ||(\cI - \cD)\Phi^{(n)}||_{V^N} + \cO ( \| ( \cI- \mathcal{D}) \Phi^{(n)} \|_{V^N}^2 )~ \leq ~ \chi  ||(\cI - \cD)\Phi^{(n)}||_{V^N}, $$
where $\chi <1 $ for $||(\cI - \cD)\Phi^{(n)}||_{V^N}$ small enough to neglect the quadratic term.
\feddich
\medskip

The convergence estimate \eref{mainconvexpr} for Algorithm 2 is easily derived from this: Consider \beqnr \cF_2(\Phi) = \Phi - (\cI - \cD_{\Phi})\mathcal{B}^{-1}(\cI-  \mathcal{D}_{\Phi} )  \cJ'(\Phi) \label{F2}, \eeqnr for which $\Phi^{(n+1)} = P (\cF_2(\Phi^{(n)}))$ for the iterates of Algorithm 2. 
Differentiation of $\cF_2$ at $\bar{\Psi}$ chosen as before gives
\beq \cF '_2(\bar{\Psi}) \Delta \Psi &=&  \cI - (\cI - \cD)\cB^{-1} (\cI - \cD) \cL^{(2)}(\bar{\Psi}, \Lambda) \Delta \Psi + \cO ( \| ( \cI- \mathcal{D}) \Phi^{(n)} \|_{V^N}^2 ),\eeq (note that derivation of the projector $D_{\bar{\Psi}}$ on the left hand side with respect to $\bar{\Psi}$ results in a zero term), so that the same reasoning as above gives
\beq&&  \|( \cI- \mathcal{D}) \Phi^{(n+1)}   \|_{V^N}\\ & \leq& || ( \cI- \mathcal{D})\big (\cI - \hat{\cB}^{-1} ( \cI- \mathcal{D})\cL^{(2, \Psi)}(\bar{\Psi}, \Lambda)\big) (\cI - \cD)\Psi ||_{V^N} + \cO ( \| ( \cI- \mathcal{D}) \Phi^{(n)} \|_{V^N}^2\\
&\leq & \chi  ||(\cI - \cD)\Phi^{(n)}||_{V^N},
\eeq
with $\chi < 1$ for $\Phi^{(n)}$ close enough to $\Psi.$
\feddich

To prove the convergence of the exponential parametrisation (Algorithm 3) defined by 
$$ \Phi^{(n+1)} := \exp\left(- \alpha \hat{\XX}\right) (\Phi^{(n)}),$$  it is enough to notice, cf. the remarks after Lemma \ref{geodesics}, that we follow a geodesic path in direction
$(I-\cD_{\Phi^{(n)}}) \cB^{-1}(A_{\Phi^{(n)}}\Phi^{(n)} - \Phi^{(n)}\mathbf{\Lambda}^{(n)})$,
which is equal to the descent direction of Algorithm 2. Due to the definition of the tangent manifold, 
 $\Phi^{(n+1)}$ again differs from $\cF_2(\Phi^{(n)})$ (defined by \eref{F2}) only by an asymptotically neglectable quadratic error term.
\begin{flushright}{$\Box$}\end{flushright}

\subsection{Quadratic convergence of the energy}

For the Rayleigh quotient $ R ( \phi^{(n)})$, i.e. for the
simplified problem and $N=1 $, it is known that $ R ( \phi^{(n)}) - R
( \psi) \lesssim \| \psi -  \phi^{(n)} \|_V^2 $. To end this section, we will 
show that this property holds also for the computed
energies, provided that the constraints are satisfied exactly and
the functional is sufficiently often differentiable. The latter is only
known for Hartree-Fock and the simplified problem. Since the exchange correlation potential
is not known exactly, this question remains open in general for the
density functional theory.

\begin{theorem}\label{thm:quad-conv}
Provided that $\cJ $ is two times differentiable on a neighborhood $U_{\delta}(\Psi)\subseteq V^N$ of the minimizer $\Psi$, and that for fixed $\Phi \in U_{\delta}(\Psi)$, $\cJ''$ is continuous on $ \{ t \Psi + (1 -t) \Phi | t \in [0, 1] \}$, the error in the energy depends quadratically on the approximation error of the minimizer $\Psi$, i.e. 
\begin{equation}\label{eq:quad-conv}
    \mathcal{J} ( \Phi ) - \mathcal{J} ( \Psi )~~ \lesssim~~
    \| ( I - \mathcal{D}_{\Psi} ) \Phi^{(n)} \|_{V^N}^2  \ .
\end{equation}
\end{theorem}

\begin{proof}
Let us choose a representant of the solution $\Psi$ according to Lemma \ref{lem:projdiffer}.
Abbreviating $e =  \Phi - \Psi $, we can use  $\cJ'(\Psi)((\cI-\cD)\Phi) = 0$ to find that $$\cJ'(\Psi)(e) ~~=~~ \cJ'(\Psi)((\cI-\cD)\Phi) ~+~ \cO(||(\cI-\cD)\Phi||^2)~~ =~~ \cO(||(\cI-\cD)\Phi||^2)$$ so that 
\beq \cJ(\Phi) -  \cJ ( \Psi ) ~&=&~ \int\limits_0^1 \cJ'(\Psi + s e)(e) ds ~+~ \frac{1}{2} \cJ' (\Phi ) (e) \\ && -~ \frac{1}{2} (\cJ' (\Psi ) (e) ~+~  \cJ' (\Phi) (e))~ +~ \cO(||(\cI-\cD)\Phi||^2).\eeq
By integration by parts, 
$$ \frac 12 ( f (0)  +  f(1) ) ~~=~~  \int\limits_0^1 f(t) dt~ +~  \int\limits_0^1 (s-\frac 12)f'(s)ds,$$
so that 
$$  \mathcal{J} ( \Phi ) - \mathcal{J} ( \Psi ) ~=~ \frac 12 \langle \langle \cJ'(\Phi)  , \Phi - \Psi \rangle \rangle ~ -~ \int\limits_0^1  (s-\frac 12) \cJ'' ( \Phi+se)(e,e)ds~+~ \cO(||(\cI-\cD)\Phi||^2).
$$ 
For estimation of the first term on the right hand side, recall from \eref{resifirst} that
$$||(\cI-\cD )\cJ'(\Phi)||_{V^N} ~~\lesssim~~ \| ( I -
\mathcal{D} ) \Phi \|_{V^N},$$
and therefore
\beq
\frac 12 \langle \langle \cJ'(\Phi)  ,  \Phi - \Psi \rangle \rangle ~&=&~ \frac 12 \langle \langle ( \cI - \cD) \cJ'(\Phi)  , ( \cI - \cD) \Phi  \rangle \rangle ~+~ \cO(\| ( I -
\mathcal{D})\Phi||^2)\\ ~&=&~ \cO(|| ( I -
\mathcal{D})\Phi||^2) , \eeq while for the second term, $| \int\limits_0^1  (s-\frac 12) \cJ'' ( \Phi+se)(e,e)ds |~=~ \cO(||e||^2) = \cO(|| ( I - \mathcal{D})\Phi||^2) $ follows from the continuity of $\cJ''$ and, again, the usage of Lemma \ref{lem:projdiffer}. 
\end{proof}

\section{Further Comments and Conclusions}\label{commentsec}

Before we conclude this article with numerical examples, we would like to make some comments about the complexity of the
numerical schemes when applied to the problems of section \ref{schrosec}, and about the potentialities for accelerating convergence of the iteration scheme.

{\bf Complexity:} Concerning disk storage, the task is to compute $N$ functions $ \psi \in V_h$, so $\mathcal{O}( N
\mbox{dim} V_h ) $ memory is needed to store the orbital functions, while storage of the
discretization of the Fock operator $ A $ requires at most $\mathcal{O}( (\mbox{dim}
V_h)^2 ) $ in the general  and worst case, but only $ \mathcal{O}(
\mbox{dim} V_h ) $ for sparse discretizations. 
Regarding computational demands, the non-zero entries of a sparse discretization of $A$ are of $\cO(\dimens V_h)$, so that the complexity of the application of $A$ depends linearily on $
\mbox{dim} V_h $. The computation of $
\langle A \hat{\phi}^{(n+1)}_i,\hat{\phi}^{(n+1)}_j \rangle $, and $
\langle  \hat{\phi}^{(n+1)}_i,\hat{\phi}^{(n+1)}_j \rangle $ needs $
\mathcal{O}( N^2   (\mbox{dim} V_h) ) $ operations in the case of sparse
discretizations (and $ \mathcal{O}( N^2  (\mbox{dim} V_h)^2 ) $ in
the worst case). The orthogonalization procedure, i.e. the
projection onto the Stiefel manifold usually has a complexity $
\mathcal{O} (N^2\mbox{dim} V_h)$. 
To relate the above complexities to the size $N$ of the electronic system, 
it is also interesting to discuss how large $\dimens V_{h,min}$ has to be chosen for a given size $N$.
To this end, we might fix a given maximal error $e$ per atom or electron (usually requested to be smaller than the intrinsic modeling error of DFT or HF models) and determine the minimal ansatz space dimension $\dimens V_{h,min}(N)$ that keeps the numerical error under that error $e$. 
If we then consider the scaling of $\dimens V_{h,min}$ with respect to the
size of the system $N$, it turns out that $\mbox{dim} V_{h,min}(N) =
\mathcal{O} (N),$ where the constant in front of $N$ is extremely
large for systematic basis functions and surprisingly small for
Gaussian type basis functions. 
Therefore, the natural scaling of the orbital based
DFT and/or HF computations with respect to the size $N$ of the underlying system gives an overall complexity of
$ \mathcal{O} (N^3)$ (or even $ \mathcal{O}(N^4)$ for non-sparse
discretizations).\\
This can be improved if the discretization of the individual
orbitals $ \phi_i^{(n)} $ requires substantially less than $
\mbox{dim} V_h $ DOFs. In an optimal case, one may archieve $ \mathcal{O} (1) $ for a fixed accuracy per atom; this is
for example the case if the diameter of support of $ \phi^{(n)}_i$ is of 
$\mathcal{O} (1) $, i.e. the support is local. In this case, the total
complexity scales only linearly with respect to $N$. Usually, the
eigenfunctions $ \breve{\psi}_i $  have global support. For
insulating materials, though, there exists a representation $ \Psi_{loc} $
such that $ [\Psi_{loc} ] = [ \breve{\Psi} ] \in \mathcal{G} $ and $ |\psi_{loc,i} ( x )| \lesssim e^{-\alpha | x-x_i | } $, $
\alpha
> >0 $ sufficiently large. These representations are called maximally localized or
Wannier orbitals. Linear scaling $ \mathcal{O} (N)$ can be achieved
if, during the iteration, the
representant $ \Phi^{(n)}_{loc} $ in the Grassmann manifold is selected
 and approximated in a way that the
diameter of support is of $\mathcal{O} (1) $.
This is the strategy pursued in Big DFT to achieve linear scaling, \cite{goedeckerhandb, goedeckerlinsc}.
We defer the further details to a forthcoming paper. A
related  approach, computing localized orbitals in an alternative way was proposed by \cite{hager} and 
exhibits extremely impressing results.

{\bf Convergence and Acceleration:} In the present paper we have considered
linear convergence of a preconditioned gradient algorithm. For the
simplified model, this convergence is guaranteed by the spectral gap
condition, in physics referred as the HOMO-LUMO gap (i.e. highest
occupied molecular orbital-lowest unoccupied molecular orbital gap). For the Hartree-Fock model,
this condition is replaced by the coercivity condition \ref{lagrangeelli}. The
same condition applies to models in density functional theory,
provided the Kohn-Sham energy functionals are sufficiently often
differentiable. Let us mention that a verification of this
conditions will answer important open problems in Hartree-Fock
theory, like uniqueness etc. The performance of the algorithm may
be improved by an optimal line search, replacing $ \mathcal{B}$ by an optimal 
$\alpha_n \mathcal{B}$. Except for the simplified problem, where an
optimal line search performed like in the Jacobi-Davidson
algorithm as a particular simple subspace acceleration, optimal line
search is rather expensive though and not used in practice.\\
Since the present preconditioned steepest decent algorithm is
gradient directed, a line search based on the Armijo rule will
guarantee convergence in principle, even without a coercivity
condition \cite{armijo, geiger}. \\
In practice, convergence is improved  by subspace acceleration
techniques, storing iterates $ \Phi^{(n-k)}, \ldots, \Phi^{(n)},
\widehat{\Phi}^{(n+1)}$ and compute $ \Phi^{(n+1)} $ from an
appropriately chosen linear  combination of them. Most prominent
examples are the DIIS  \cite{diis} and conjugate gradient
\cite{cg, arias} algorithm.The DIIS algorithm is implemented in the EU NEST
project BigDFT, and
frequently used in other quantum chemistry codes. Without going into
detailed descriptions of those methods and further investigations,
let us point out that the analysis in this paper provides the
convergence of the worst case scenario.
Second order methods, in particular Newton methods  have been
proposed in literature \cite{maday}, but since these require
the solution of a linear system of size  $ N \mbox{dim} V_h \times N
\mbox{dim} V_h $, they are to be avoided.

\section{Numerical examples} 

 \begin{figure}
 \begin{center}
 \includegraphics[width=0.4\textwidth]{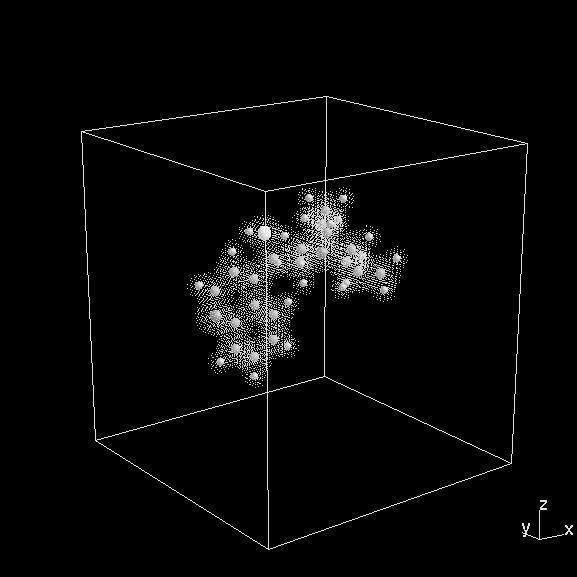}
\hspace{0.05\linewidth}
 \includegraphics[width=0.4\textwidth]{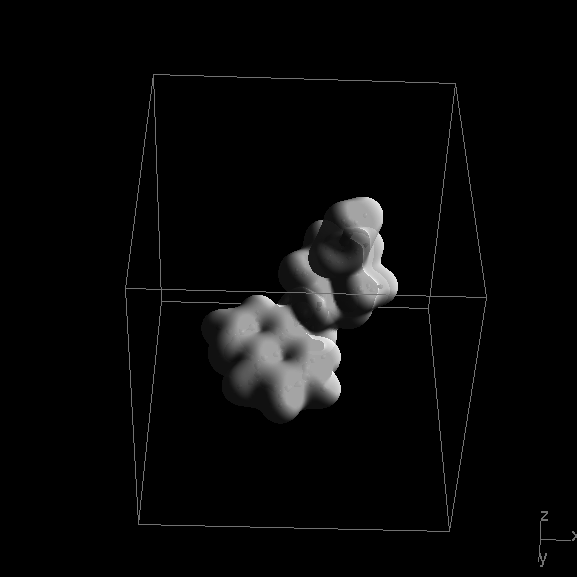}
 \caption{Atomic geometry and electronic structure of cinchonidine\label{geometry}}
 \end{center}
 \end{figure}

The proposed direct minimization algorithm 1 is realized in the recent density
functional code bigDFT \cite{bigdft}, which is implemented in the open source ABINIT package, a common project of the Universit\'{e} Catholique de Louvain, Corning Incorporated, and other contributors \cite{abinit, abinit2, abinit3, bigdftpap}.
It relies on an efficient Fast Fourier Transform algorithm \cite{goedefft} for the conversion of wavefunctions between real and reciprocal space, together with a DIIS subspace acceleration.
 We demonstrate the convergence for the simple
molecule cinchonidine ($C_{19}H_{22}N_2O $) of moderate size $N = 55$
for a given geometry of the nuclei displayed in figure \ref{geometry}. 
Despite
the fact that the underlying assumptions in the present paper cannot
 be verified rigorously, the proposed convergence behavior is
observed by all benchmark computations. The  algorithm is
experienced to be quite robust also if the HOMO-LUMO gap is
relatively small. \\

\begin{figure}
\begin{center}
\epsfig{file=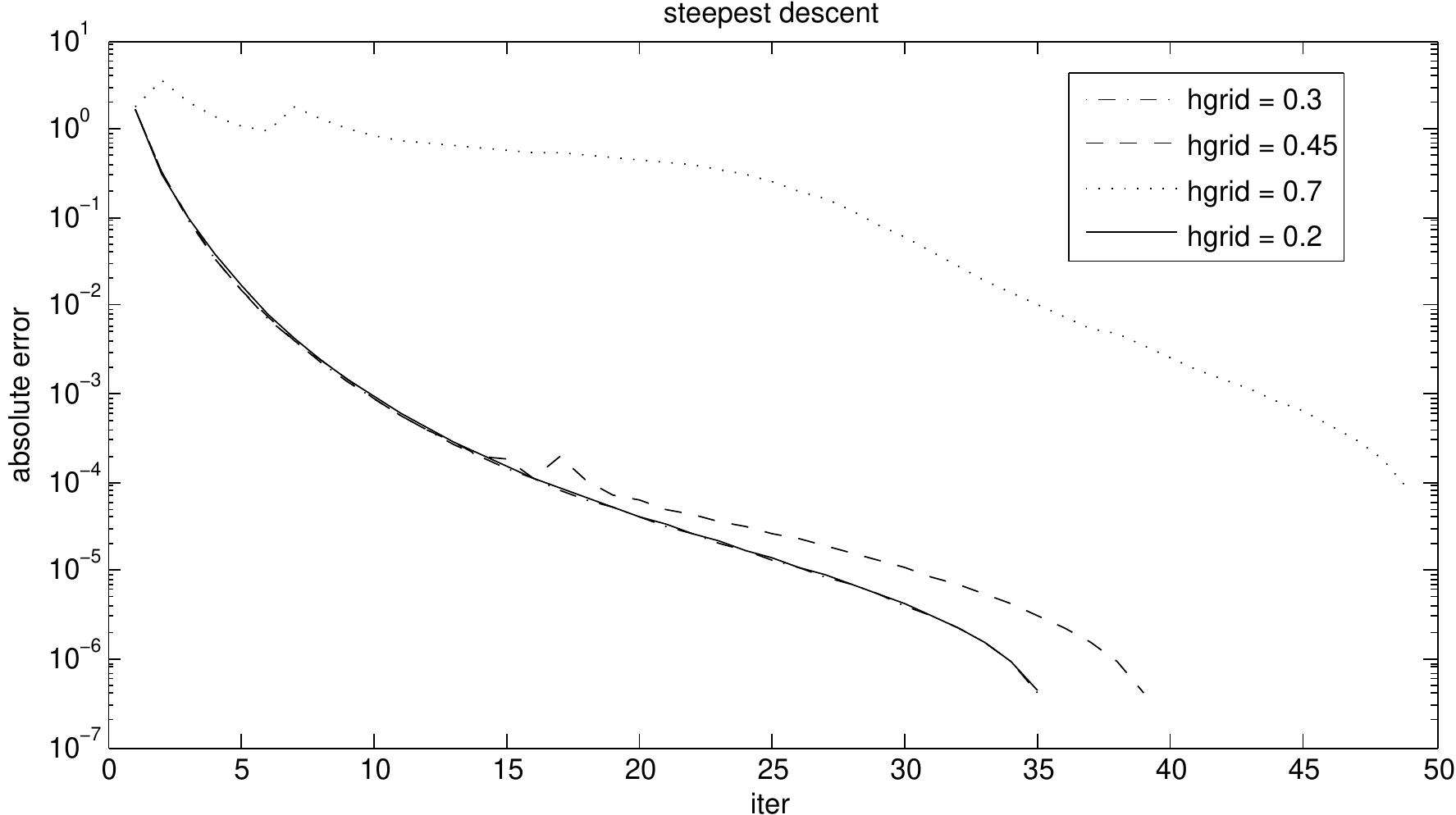, angle=0, width=0.45\linewidth, height=0.3\linewidth}
\hspace{0.05\linewidth}
\epsfig{file=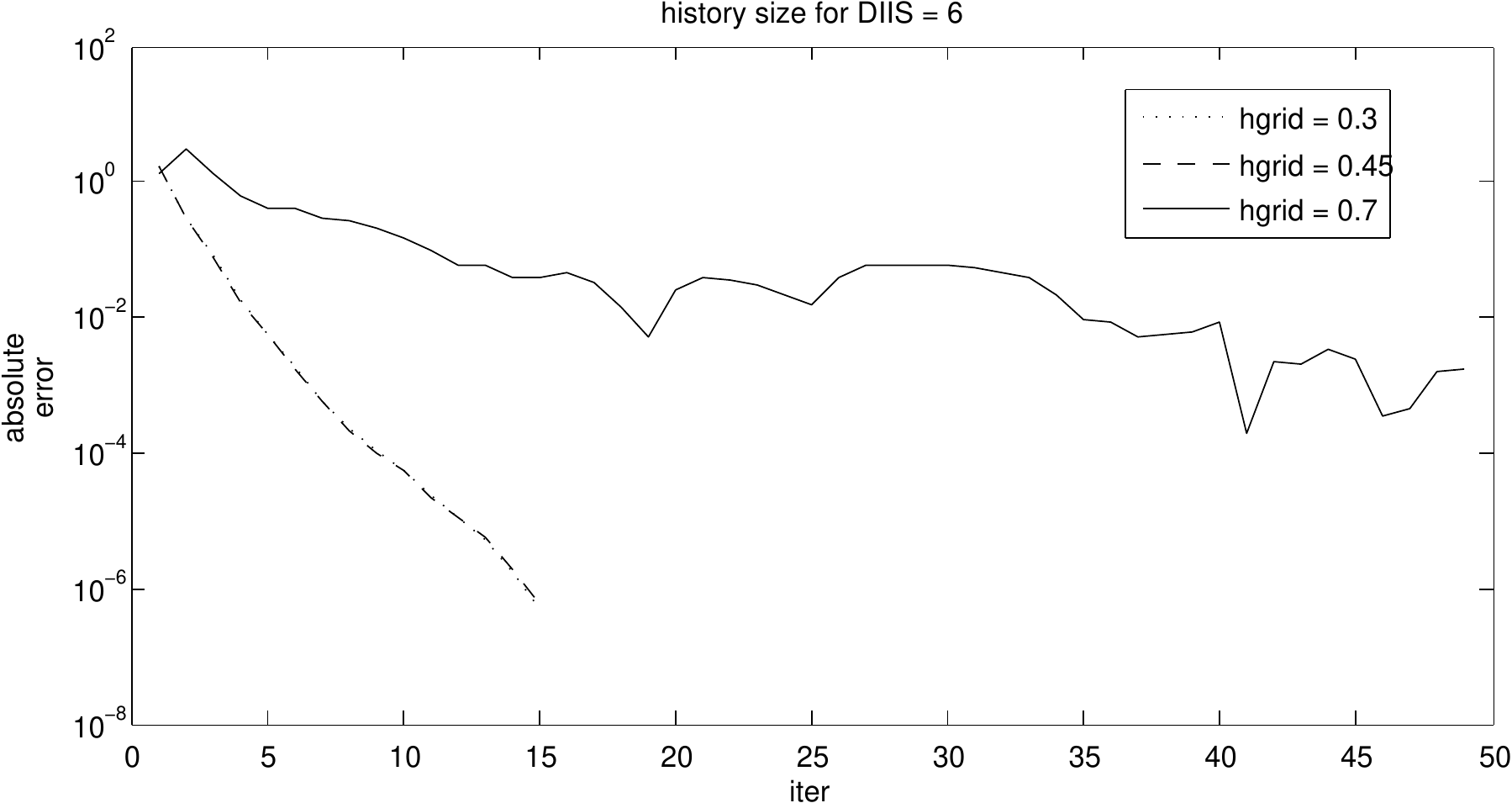, angle=0, width=0.45\linewidth, height=0.3\linewidth}
\caption{\label{convhist}Convergence history for the direct minimization scheme (left)  and with DIIS acceleration (right) for different mesh sizes.}
\bigskip

\epsfig{file=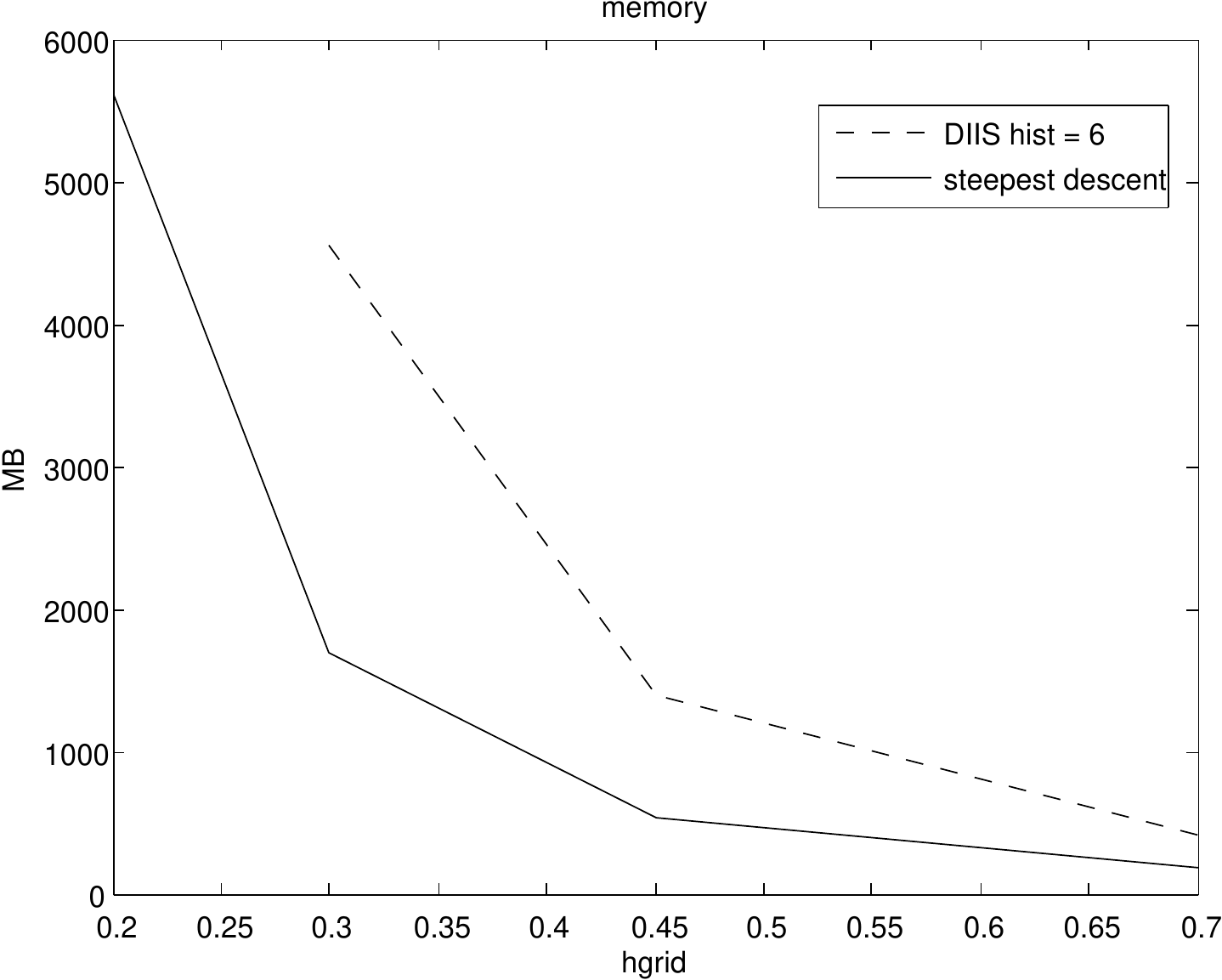, angle=0, width=0.45\linewidth, height=0.3\linewidth}
\hspace{0.05\linewidth}
\epsfig{file=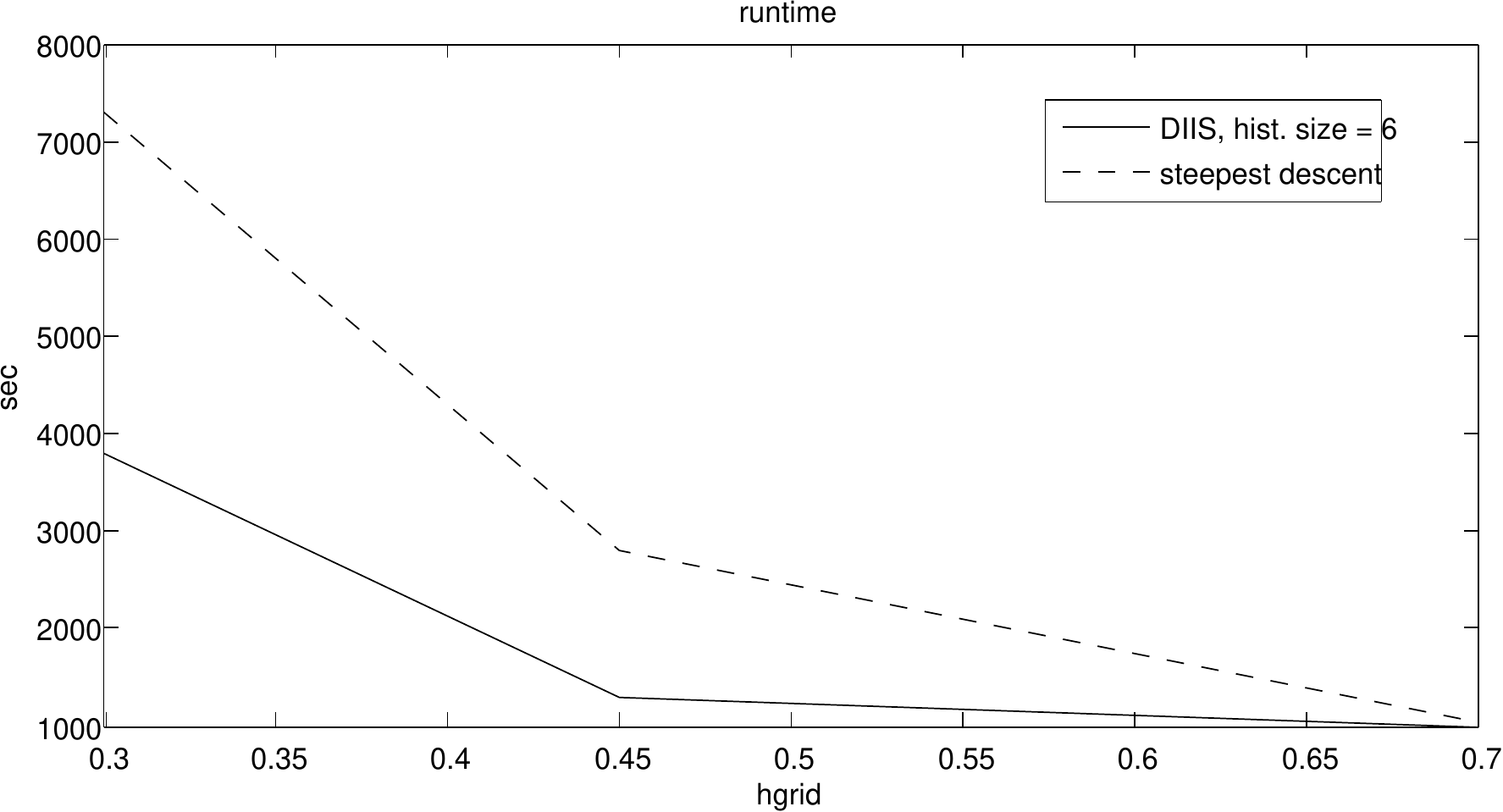, angle=0, width=0.45\linewidth, height=0.3\linewidth}

\caption{Memory requirements (left) and computing time (right) for direct minimization algorithm with and without DIIS acceleration.}
\end{center}
\end{figure}
For our computations, we have used a simple LDA (local density approximation) model
proposed by \cite{goedecker-umrigar} and norm-conserving non-local
pseudopotentials \cite{goedecker-hutter}. The  orbital functions $
\psi_i $ are approximated by Daubechies orthogonal wavelets with $8$
vanishing moments based on an approximate Galerkin discretization \cite{goedeckermagic}.
For updating the nonlinear potential, the electron density is
approximated by interpolating scaling functions (of order 16). The
discretization error can be controlled by an underlying basic mesh
size $ h_{grid} $.  

In figure \ref{convhist}, we demonstrate the convergence 
of the present algorithm for $ 4 $ different choices of mesh sizes, where the error is given in the energy norm of the discrete functions. The initial guess for the orbitals is given by the atomic solutions. Except in case of non-sufficient resolution ($h_{grid} = 0.7 $), where we
obtain a completely wrong result, convergence is observed. If the discretisation is sufficiently good, we do not observe much difference in the convergence history for different mesh sizes. Since the convergence speed depends on the actual solution, it is only possible to observe that
the convergence is bounded by a linear rate.

The number
of iterations is relatively moderate bearing in mind that one
iteration step only requires matrix-vector multiplications with the Fock
operator and not a corresponding solution of linear equations. 
The DIIS implemented in BigDFT accelerates the iteration by almost
halving the  number of iterations and the total computing time at the expense
of additional storage capacities, see also figure \ref{convhist}. 
Further benchmark computations have already been performed and will be reported in different publications
by the groups involved in the implementation of BigDFT.

\end{document}